\documentclass[12pt]{article}

\usepackage{}

\newtheorem{theorem}{Theorem}[section]
\newtheorem{lemma}[theorem]{Lemma}

\newtheorem{definition}[theorem]{Definition}

\begin{document}

\title{\vspace*{-1.5cm}
Strong Convergence of Wong-Zakai Approximations of Reflected SDEs in A Multidimensional General Domain }

\author{Tusheng Zhang$^{1}$}
\footnotetext[1]{\  School of Mathematics, University of Manchester,
Oxford Road, Manchester M13 9PL, England, U.K. Email: 
tusheng.zhang@manchester.ac.uk}
\maketitle

\begin{abstract}
 In this paper, we obtained the strong convergence of Wong-Zakai approximations of reflected SDEs in a general multidimensional domain  giving an affirmative answer to the question posed in [ES].
\end{abstract}

\section{Introduction}
 Let $D$ be a bounded domain in $R^d$. Consider the reflected stochastic differential equation (SDE):
 \begin{eqnarray} \label{1.1}
 \left\{\begin{array}{l}dX(t)=\sigma (X(t))\circ dW(t)+b(X(t))dt+dL(t),\\
   X(0)=x_0, \quad X(t)\in \bar{D}, t\geq 0,\\
   |L|(t)=\int_0^tI_{\partial D}(X(s))d|L|(s),
   \end{array}
   \right.
 \end{eqnarray}
 where $W(t), t\geq 0$ is a $m$-dimensional Brownian motion, $|L|(t)$ stands for the total variation of $L$ on the interval $[0,t]$, $\circ$ indicates a
 Stratonovich integral.
 \vskip 0.3cm
 There is a big amount of literature devoted to the study of reflected SDEs. Let us mention a few of them.
 Reflected SDEs in a convex domain was first studied by H. Tanaka in \cite{T}. Existence and uniqueness of solutions of reflected SDEs in  general domains were established by Lions and Sznitman in \cite{LS} and  Saisho in \cite{S}. Existence and uniqueness of solutions of reflected SDEs under more general coefficients than the usual Lipschitz conditions were considered in \cite{MR}.

 \vskip 0.3cm
 The purpose of this paper is to study Wong-Zakai type approximations of above reflected SDEs. Let $W^n$ be the $n-$dyadic piecewise linear interpolation of $W$ and  $X^n$  the solution of the following reflected random ordinary
 differential equation:
\begin{eqnarray} \label{1.2}
 \left\{\begin{array}{l}\dot{X}^n(t)=\sigma (X^n(t))\dot{W}^n(t)+b(X^n(t))dt+\dot{L}^n(t),\\
 X^n(0)=x_0, \quad X^n(t)\in \bar{D}, t\geq 0,\\
   |L^n|(t)=\int_0^tI_{\partial D}(X^n(s))d|L^n|(s).
   \end{array}
   \right.
 \end{eqnarray}

We are concerned with the strong  convergence of  $X^n$ to the solution $X$. Strong convergence of Wong-Zakai approximations to stochastic differential equations is well known, see e.g. \cite{IW}. However, the convergence of  Wong-Zakai approximations to stochastic differential equations with reflection (especially in higher  dimensions) is quite tricky because of the constraints on the solution and the appearance of the boundary local time. As far as we are aware of, there are two main papers related to this question. In \cite{P}, Petterson established a Wong-Zakai approximations for SDEs with reflection under the assumption that the domain is convex. The convexity is too rigid sometimes for applications. In \cite{ES}, Evans and Stroock considered Wong-Zakai approximations for reflected SDEs in general domains (as in \cite{LS}) and proved that
$X^n$ converges weakly (in law ) to the solution $X$. In the same paper, the authors also posed the question of whether the strong convergence holds. For some of the interesting applications, we refer the reader to \cite{ES}.
\vskip 0.3cm
The purpose of this paper is to establish the strong convergence ( the $L^p$ convergence in $C([0,T], \bar{D})$ of the Wong-Zakai approximations for reflected SDEs in multidimensional general domains, hence giving an affirmative answer to the question in \cite{ES}.

\vskip 0.3cm
The paper is organized as follows. In Section 2, we recall the framework and formulate the main result. The rest of the paper ( Section 3 ) is entirely  devoted to the proof of the theorem.

\section{Framework and the main result}
\setcounter{equation}{0}
Let $D\subset R^d$ be a bounded domain with boundary $\partial D$. For $x\in \partial D$, let $\nu(x)\subset S^{d-1}$ denote a nonempty collection of reflecting directions. Throughout this paper, as in [LS], [ES], we impose the following conditions on the domain.
\begin{itemize}
\item[\bf{D.1}] $\nu(x)\not=\emptyset$ for every $x\in \partial D$ and there exist a constant $C_0\geq 0$ such that
$$(x^{\prime}-x)\cdot \nu +C_0|x-x^{\prime}|^2\geq 0\quad\quad \mbox{for all}\quad x^{\prime}\in D, x\in \partial D\quad \mbox{and}\quad  \nu\in \nu(x).$$
\item[\bf{D.2}] There exists a function $\phi\in C^2(R^d; R)$ and $\alpha >0$ such that
$$\nabla \phi(x)\cdot \nu\geq \alpha \quad\quad \mbox{for all}\quad x\in \partial D\quad \mbox{and}\quad  \nu\in \nu(x).$$
\item[\bf{D.3}] There exist $n\geq 1$, $\lambda>0$, $K>0$, $a_1,a_2,...,a_n\in S^{d-1}$, and $x_1,x_2,...,x_n\in \partial D$ such that  $\partial D\subset \cup_{i=1}^n B(x_i, K)$ and $x\in \partial D\cap B(x_i, 2K)\Longrightarrow \nu\cdot a_i\geq \lambda$ for all $\nu\in \nu(x)$.
\end{itemize}
\noindent{\bf Convention}; Throughout this paper,  any function $G$ defined on the positive half line $[0,\infty)$ automatically extends to a function on the whole line by setting $G(s)=G(s\vee 0)$ when necessary.
\vskip 0.3cm

Let $W(t)=(W_1(t),W_2(t),...,W_m(t)), t\geq 0$ be a $m$-dimensional Brownian motion on a completed filtered probability space $(\Omega, {\cal F}, {\cal F}_t, P)$. Suppose $\sigma=(\sigma_{i,j}) \in C^1(\bar{D}; R^d\otimes R^m)$ such that the derivative $\sigma^{\prime}$ is Lipschitz continuous and that $b: \bar{D}\rightarrow R^d$ is Lipschitz
continuous.
\vskip 0.3cm
For $n\in N$ and $s\in [\frac{k}{2^n}, \frac{k+1}{2^n})$, set $s_n^-=(\frac{k-1}{2^n})\vee 0$ and $s_n=\frac{k}{2^n}$. Let $W^n$ be the linear interpolation of $W$ defined by
\begin{equation} \label{2.1}
 W^n(t)=W(\frac{k-1}{2^n})+2^n(t-\frac{k}{2^n})(W(\frac{k}{2^n})-W(\frac{k-1}{2^n}))
 \end{equation}
 for  $t\in [\frac{k}{2^n}, \frac{k+1}{2^n}), \quad k=0,1,2,...$
 Note that the above convention applies here.
 Let $\sigma\sigma^{\prime}: \bar{D}\rightarrow R^d$ be defined as
 \begin{equation} \label{2.2}
 (\sigma\sigma^{\prime}(y))_i=\sum_{j=1}^m\sum_{k=1}^d \frac{\partial \sigma_{i,j}(y)}{\partial y_k}\sigma_{k,j}(y).
 \end{equation}
 With this notation, equation (\ref{1.1}) becomes
 \begin{equation} \label{2.3}
 X(t)=x_0+\int_0^t\sigma (X(s))dW(s)+ \frac{1}{2}\int_0^t\sigma\sigma^{\prime}(X(s))ds +\int_0^tb(X(s))ds+L(t)
 \end{equation}
 \begin{definition}
 We say that $(X,L)$ is a solution to the reflected SDE (\ref{2.3}) if $(X,L)$ is a $\bar{D}\times R^d$-valued, adapted continuous
 process such that

 (i) $L(t), t\geq 0$ is of bounded variation on any finite sub-interval of $[0,\infty )$,

 (ii) for $t\geq 0$,
 $$ X(t)=x_0+\int_0^t\sigma (X(s))dW(s)+ \frac{1}{2}\int_0^t\sigma\sigma^{\prime}(X(s))ds +\int_0^tb(X(s))ds+L(t)$$
 almost surely,

 (iii)

 $$|L|(t)=\int_0^tI_{\partial D}(X(s))d|L|(s), \quad \quad \quad L(t)=\int_0^t\nu(X(s))d|L|(s),$$
 where $|L|(t)$ stands for the total variation of $L$ on the interval $[0,t]$, the last equality means that $\frac{DL(t)}{d|L|(t)}\in \nu(X(t))$.
 \end{definition}
 The solution $(X^n, L^n)$ to the  reflected random ordinary
 differential equation (\ref{1.2}) is defined accordingly.
 \vskip 0.3cm
 Under the above assumptions, the existence and uniqueness of $X^n, X$ are well known now, see, for example, \cite{LS}.
Here is the main result.
\begin{theorem}
Let $X^n,X$ be the solutions to reflected stochastic equations (\ref{1.1}) and (\ref{1.2}). It holds that for any $p>0$ and $T>0$,
\begin{equation}\label{2.4}
\lim_{n\rightarrow \infty}E[\sup_{0\leq t\leq T}|X^n(t)-X(t)|^p]=0.
\end{equation}
\end{theorem}

In next section, $C$ will denote a generic constant which is usually different from line to line.

\section{The
roof of the main result}
\setcounter{equation}{0}
The rest of the paper is devoted to the proof of Theorem 2.2. First of all we recall the following estimate from \cite{ES}.
\begin{lemma}
Let $p\geq 2$, $T>0$. Then there exists a constant $C_1(T,p)$ independent of $n$ such that
\begin{equation}\label{3.01}
E[|X^n(t)-X^n(s)|^p]\leq C_1(T,p) |t-s|^{\frac{p}{2}},
\end{equation}
for $0\leq s, t\leq T$.
\end{lemma}
Applying the proof in \cite{ES} to $X(t), t\geq 0$ one also has
\begin{lemma}
Let $p\geq 2$, $T>0$. Then there exists a constant $C_2(T,p)$ such that
\begin{equation}\label{3.02}
E[|X(t)-X(s)|^p]\leq C_2(T,p) |t-s|^{\frac{p}{2}},
\end{equation}
for $0\leq s, t\leq T$.
\end{lemma}
Due to (\ref{3.01}), (\ref{3.02}) above, to prove Theorem 2.2, it is sufficient to show that for any fixed $t>0$
\begin{equation}\label{3.03}
\lim_{n\rightarrow \infty}E[|X^n(t)-X(t)|^2]=0.
\end{equation}
Indeed, it follows from (\ref{3.01}), (\ref{3.02}) and Garsia, Rodemich and Rumsey's lemma (See Theorem 1.1 in \cite{W}) that for a fixed positive number $\alpha_0<\frac{1}{2}$, there exist random variables $K_n(\omega), K(\omega)$ such that
\begin{equation}\label{3.06}
|X^n(t)-X^n(s)|\leq K_n(\omega) |t-s|^{\alpha_0}, \quad s, t\in [0,T],
\end{equation}
and
\begin{equation}\label{3.07}
|X(t)-X(s)|\leq K(\omega)|t-s|^{\alpha_0},\quad s, t\in [0,T].
\end{equation}
Furthermore, because the constant $C_1(T,p)$ in (\ref{3.01}) is independent of $n$, $K_n, K$ can be chosen to satisfy
\begin{equation}\label{3.08}
\sup_nE[K_n^p]<\infty, \quad E[K^p]<\infty,
\end{equation}
for any $p>0$. Since $X^n, X$ live on the bounded domain $\bar{D}$, to show (\ref{2.4}) it is sufficient to prove that $X^n$ converges to $X$ in probability. For $\varepsilon>0$, choose $t_i\in [0,T], i=1,...,N_{\varepsilon}$ such that $[0,T]\subset \cup_iB(t_i, \varepsilon)$. Given $\delta>0$, for $M>0$, we have
\begin{eqnarray}\label{3.09}
&&P(\sup_{0\leq t\leq T}|X^n(t)-X(t)|>\delta )\nonumber\\
&\leq& P(\sup_i\sup_{t\in B(t_i,\varepsilon)}|X^n(t)-X^n(t_i)|>\frac{\delta}{3}, |K_n|\leq M)+P(|K_n|> M)\nonumber\\
&+&P(\sup_i\sup_{t\in B(t_i,\varepsilon)}|X(t)-X(t_i)|>\frac{\delta}{3}, |K|\leq M)+P(|K|> M)\nonumber\\
&+& \sum_{i=1}^{N_{\varepsilon}}P(|X^n(t_i)-X(t_i)|>\frac{\delta}{3}).
\end{eqnarray}
Now, for any given $\eta>0$, by (\ref{3.08}) we first choose $M$ sufficiently large so that $P(|K_n|> M)\leq \frac{\eta}{4}$, $P(|K|> M)\leq \frac{\eta}{4}$ for all $n$. For such a constant $M$, because of (\ref{3.06}) and (\ref{3.07}) we can select $\varepsilon>0$ sufficiently small so that
$$P(\sup_i\sup_{t\in B(t_i,\varepsilon)}|X^n(t)-X^n(t_i)|>\frac{\delta}{3}, |K_n|\leq M)=0,$$
and
$$P(\sup_i\sup_{t\in B(t_i,\varepsilon)}|X(t)-X(t_i)|>\frac{\delta}{3}, |K|\leq M)=0.$$
When $\varepsilon$ is fixed, it follows from (\ref{3.03}) that there exists $N>0$ such that for $n\geq N$,
$$\sum_{i=1}^{N_{\varepsilon}}P(|X^n(t_i)-X(t_i)|>\frac{\delta}{3})<\frac{\eta}{4}$$
Putting the above arguments together we prove that $X^n$ converges to $X$ in probability.
\vskip 0.3cm
So we remain to prove (\ref{3.03}).
Again because of (\ref{3.01}), (\ref{3.02}) we may assume that $t$ is a dyadic number, i.e.,  $t=\frac{k_0}{2^{n_0}}$ for some positive integers $k_0$, $n_0$ and we may also assume $n\geq n_0$.
\vskip 0.3cm
Let $f(y_1,y_2,y_3)=exp(r(y_1+y_2))y_3$. Recall $\phi$ is the function specified in (D.2). To simplify the exposure, we introduce  the following notation:
$$y_1(t):=\phi(X(t)), y_2^n(t):=\phi(X^n(t)), y_3^n(t):=|X^n(t)-X(t)|^2.$$
$$f_n(t):=f(y_1(t),y_2^n(t),y_3^n(t)), g_n(t):=exp(ry_1(t)+ry_2^n(t)).$$
Since $X^n, X$ take values in the bounded domain $\bar{D}$, we have
\begin{equation}\label{3.04}
c_1 |X^n(t)-X(t)|^2 \leq  f_n(t)\leq c_2|X^n(t)-X(t)|^2,
\end{equation}
where $c_1, c_2$ are positive constants independent of $n$.
Thus the proof of (\ref{3.03}) reduces to  show
\begin{equation}\label{3.05}
\lim_{n\rightarrow \infty}E[f_n(t)]=0.
\end{equation}
\vskip 0.3cm
By Ito's formula, we have
\begin{eqnarray}\label{3.1}
&& f_n(t)\nonumber\\
&=& r\int_0^tf_n(s)<\nabla\phi(X(s)), \sigma(X(s))dW(s)>+r\int_0^tf_n(s)<\nabla\phi(X(s)), b(X(s)>ds\nonumber\\
&+& \frac{1}{2}r\int_0^tf_n(s) tr(\phi^{\prime\prime}(\sigma\sigma^*)(X(s)))ds+\frac{1}{2}r\int_0^tf_n(s) <\nabla \phi(X(s),\sigma\sigma^{\prime}(X(s))> ds\nonumber\\
&+& r\int_0^tf_n(s)<\nabla\phi(X(s)), \nu(X(s))>d|L|(s)+r\int_0^tf_n(s)<\nabla\phi(X^n(s)), \sigma(X^n(s))dW^n(s)>\nonumber\\
&+&r\int_0^tf_n(s)<\nabla\phi(X^n(s)), b(X^n(s))>ds+r\int_0^tf_n(s)<\nabla\phi(X^n(s)), \nu(X^n(s))>d|L^n|(s)\nonumber\\
&+&2\int_0^tg_n(s)<X^n(s)-X(s), \sigma(X^n(s))dW^n(s)>\nonumber\\
&-&2\int_0^tg_n(s)<X^n(s)-X(s), \sigma(X(s))dW(s)>\nonumber\\
&+&2\int_0^tg_n(s)<X^n(s)-X(s), b(X^n(s))-b(X(s))>ds\nonumber\\
&-&\int_0^tg_n(s)<X^n(s)-X(s), \sigma\sigma^{\prime}(X(s))>ds\nonumber\\
&+&2\int_0^tg_n(s)<X^n(s)-X(s), \nu(X^n(s))d|L^n|(s)-\nu(X(s))d|L|(s)>\nonumber\\
&+&\int_0^tg_n(s) tr(\sigma\sigma^*(X(s)))ds+\frac{1}{2}r^2\int_0^tf_n(s)|\sigma^*\nabla\phi|^2(X(s))ds\nonumber\\
&-&2r\int_0^tg_n(s)<\sigma^*(X(s))(X^n(s)-X(s)), \sigma^*\nabla\phi(X(s))>ds.
\end{eqnarray}

\begin{eqnarray}\label{3.2}
&& g_n(t)\nonumber\\
&=& exp(2r\phi(x_0))+r\int_0^tg_n(s)<\nabla\phi(X(s)), \sigma(X(s))dW(s)>\nonumber\\
&+&r\int_0^tg_n(s)<\nabla\phi(X(s)), b(X(s))>ds\nonumber\\
&+& \frac{1}{2}r\int_0^tg_n(s) tr(\phi^{\prime\prime}(\sigma\sigma^*)(X(s)))ds+\frac{1}{2}r\int_0^tg_n(s) <\nabla \phi(X(s),\sigma\sigma^{\prime}(X(s))> ds\nonumber\\
&+& r\int_0^tg_n(s)<\nabla\phi(X(s)), \nu(X(s))>d|L|(s)\nonumber\\
&+& r\int_0^tg_n(s)<\nabla\phi(X^n(s)), \sigma(X^n(s))dW^n(s)>\nonumber\\
&+&r\int_0^tg_n(s)<\nabla\phi(X^n(s)), b(X^n(s))>ds+r\int_0^tg_n(s)<\nabla\phi(X^n(s)), \nu(X^n(s))>d|L^n|(s)\nonumber\\
&+&\frac{1}{2}r^2\int_0^tg_n(s)|\sigma^*\nabla\phi|^2(X(s))ds
\end{eqnarray}
\vskip 0.3cm
To bound $E[f_n(t)]$, the crucial step is to get proper estimates for the terms
$$ rE[\int_0^tf_n(s)<\nabla\phi(X^n(s)), \sigma(X^n(s))dW^n(s)>],$$
and
$$ rE[\int_0^tg_n(s)<X^n(s)-X(s), \sigma(X^n(s))dW^n(s)>].$$
This will be done in the following two lemmas.
\vskip 0.3cm

\begin{lemma}
It holds that
\begin{eqnarray}\label{3.3}
&& rE[\int_0^tf_n(s)<\nabla\phi(X^n(s)), \sigma(X^n(s))dW^n(s)>]\nonumber\\
&\leq & C(\frac{1}{2^n})^{\frac{1}{2}}+r^2E[\int_0^tf_n(s)<\sigma^*\nabla\phi(X(s)),\sigma^*\nabla\phi(X^n(s))>ds]\nonumber\\
&&+\frac{1}{2}r^2E[\int_0^tf_n(s)|\sigma^*\nabla\phi|^2(X^n(s))ds]\nonumber\\
&&+ r\int_0^t <g_n(s)\sigma^*(X^n(s))(X^n(s)
 -X(s)), \sigma^*\nabla\phi(X^n(s))>ds\nonumber\\
 &&+\frac{1}{2}r\int_0^tf_n(s)\sum_{i=1}^m(\sigma^*(\nabla(\sigma^*\nabla\phi)_i))_i(X^n(s)) ds\nonumber\\
 &&- 2r\int_0^t <g_n(s)\sigma^*(X(s))(X^n(s)
 -X(s)), \sigma^*\nabla\phi(X^n(s))>ds.
\end{eqnarray}
\end{lemma}
\vskip 0.3cm
\noindent{\bf Proof}.  Set
$$A=r\int_0^tf_n(s)<\nabla\phi(X^n(s)), \sigma(X^n(s))dW^n(s)>.$$
Write
\begin{eqnarray}\label{3.4}
A&=& r\int_0^tf_n(s_n^{-})<\nabla\phi(X^n(s_n^{-})), \sigma(X^n(s_n^{-}))dW^n(s)>\nonumber\\
&+&r\int_0^t(f_n(s)-f_n(s_n^{-}))<\nabla\phi(X^n(s)), \sigma(X^n(s))dW^n(s)>\nonumber\\
&+& r\int_0^tf_n(s_n^{-})<\sigma^*\nabla\phi(X^n(s))-\sigma^*\nabla\phi(X^n(s_n^{-})), dW^n(s)> \nonumber\\
&:=& A_1+A_2+A_3.
\end{eqnarray}
As a stochastic integral, it is easy to see that $E[A_1]=0$. In view of (\ref{3.1}), we further write $A_2$ as
\begin{eqnarray}\label{3.5}
&&A_2\nonumber\\
&=& r^2\int_0^t(\int_{s_n^-}^sf_n(u)<\nabla\phi(X(u)), \sigma(X(u))dW(u)>)<\nabla\phi(X^n(s)), \sigma(X^n(s))dW^n(s)> \nonumber\\
&+&r^2\int_0^t(\int_{s_n^-}^sf_n(u)<\nabla\phi(X(u)), b(X(u))du>)<\nabla\phi(X^n(s)), \sigma(X^n(s))dW^n(s)> \nonumber\\
&+&\frac{1}{2}r^2\int_0^t(\int_{s_n^-}^sf_n(u)tr(\phi^{\prime\prime}(\sigma\sigma^*)(X(u))du)<\nabla\phi(X^n(s)), \sigma(X^n(s))dW^n(s)> \nonumber\\
&+&\frac{1}{2}r^2\int_0^t(\int_{s_n^-}^sf_n(u)<\nabla \phi (X(u)), (\sigma\sigma^{\prime})(X(u))>du)<\nabla\phi(X^n(s)), \sigma(X^n(s))dW^n(s)> \nonumber\\
&+&r^2\int_0^t(\int_{s_n^-}^sf_n(u)<\nabla \phi (X(u)), \nu(X(u))>d|L|(u))<\nabla\phi(X^n(s)), \sigma(X^n(s))dW^n(s)> \nonumber\\
&+&r^2\int_0^t(\int_{s_n^-}^sf_n(u)<\nabla\phi(X^n(u)), \sigma(X^n(u))dW^n(u)> )<\nabla\phi(X^n(s)), \sigma(X^n(s))dW^n(s)> \nonumber \\
&+&r^2\int_0^t(\int_{s_n^-}^sf_n(u)<\nabla\phi(X^n(u)), b(X^n(u))du>)<\nabla\phi(X^n(s)), \sigma(X^n(s))dW^n(s)> \nonumber\\
&+&r^2\int_0^t(\int_{s_n^-}^sf_n(u)<\nabla \phi (X^n(u)), \nu(X^n(u))>d|L^n|(u))<\nabla\phi(X^n(s)), \sigma(X^n(s))dW^n(s)> \nonumber\\
&+&2r\int_0^t(\int_{s_n^-}^sg_n(u)<X^n(u)-X(u), \sigma(X^n(u))dW^n(u)>)<\nabla\phi(X^n(s)), \sigma(X^n(s))dW^n(s)> \nonumber\\
&-&2r\int_0^t(\int_{s_n^-}^sg_n(u)<X^n(u)-X(u), \sigma(X(u))dW(u)>)<\nabla\phi(X^n(s)), \sigma(X^n(s))dW^n(s)> \nonumber\\
&+&2r\int_0^t(\int_{s_n^-}^sg_n(u)<X^n(u)-X(u), b(X^n(u))-b(X(u))>du)<\nabla\phi(X^n(s)),\nonumber\\
&&\quad\quad\quad\quad  \sigma(X^n(s))dW^n(s)> \nonumber\\
&-&r\int_0^t(\int_{s_n^-}^sg_n(u)<X^n(u)-X(u), \sigma\sigma^{\prime}(X(u))>du)<\nabla\phi(X^n(s)), \sigma(X^n(s))dW^n(s)> \nonumber\\
&+&2r\int_0^t(\int_{s_n^-}^sg_n(u)<X^n(u)-X(u), \nu(X^n(u))>d|L^n|(u))\nonumber\\
&&\quad\quad \times <\nabla\phi(X^n(s)), \sigma(X^n(s))dW^n(s)> \nonumber\\
&-&2r\int_0^t(\int_{s_n^-}^sg_n(u)<X^n(u)-X(u), \nu(X(u))>d|L|(u))<\nabla\phi(X^n(s)), \sigma(X^n(s))dW^n(s)> \nonumber\\
&+&r\int_0^t(\int_{s_n^-}^sg_n(u)tr(\sigma\sigma^*(X(u)))du)<\nabla\phi(X^n(s)), \sigma(X^n(s))dW^n(s)> \nonumber\\
&+&\frac{1}{2}r^3\int_0^t(\int_{s_n^-}^sf_n(u)|\sigma^*\nabla\phi|^2(X(u)))du)<\nabla\phi(X^n(s)), \sigma(X^n(s))dW^n(s)> \nonumber\\
&-&2r^2\int_0^t(\int_{s_n^-}^sg_n(u)<\sigma^*(X(u))(X^n(u)-X(u)), \sigma^*\nabla\phi(X(u))>du)\nonumber\\
&&\quad\quad\quad \times <\nabla\phi(X^n(s)), \sigma(X^n(s))dW^n(s)> \nonumber\\
&:=& \sum_{i=1}^{17}A_{2i}
\end{eqnarray}
We will bound each of the terms. Since $\nabla\phi$, $b$, $\sigma$ are bounded on $\bar{D}$, we have
\begin{eqnarray}\label{3.6}
E[|A_{22}|]&\leq& C\int_0^t (s-s_n^-)E[|\dot{W}^n(s)|]ds\nonumber\\
&\leq& C\frac{1}{2^n}\int_0^t (2^n)^{\frac{1}{2}}ds\leq C (\frac{1}{2^n})^{\frac{1}{2}}.
\end{eqnarray}
Similarly, it holds that
\begin{eqnarray}\label{3.7}
E[|A_{2i}|]&\leq& C (\frac{1}{2^n})^{\frac{1}{2}}, \quad i=3,4,7, 11,12, 15,16,17.
\end{eqnarray}
To bound $A_{21}$, we write it as
 \begin{eqnarray}\label{3.8}
 &&A_{21}\nonumber\\
 &=&
r^2\int_0^t(\int_{s_n^-}^s[f_n(u)<\nabla\phi(X(u)), \sigma(X(u)dW(u)>-f_n(s_n^-)<\nabla\phi(X(s_n^-)),\nonumber\\
&&\quad\quad  \sigma(X(s_n^-))dW(u)>])<\nabla\phi(X^n(s)), \sigma(X^n(s))dW^n(s)> \nonumber\\
&+& r^2\int_0^t(\int_{s_n^-}^sf_n(s_n^-)<\nabla\phi(X(s_n^-)), \sigma(X(s_n^-))dW(u)>)\nonumber\\
&&\times [<\nabla\phi(X^n(s)), \sigma(X^n(s))dW^n(s)>-<\nabla\phi(X^n(s_n^-)), \sigma(X^n(s_n^-))dW^n(s)>]\nonumber\\
&+& r^2\int_0^tf_n(s_n^-)<\nabla\phi(X(s_n^-)), \sigma(X(s_n^-))(W(s)-W(s_n^-))>\nonumber\\
&&\quad\quad \times <\nabla\phi(X^n(s_n^-)), \sigma(X^n(s_n^-))dW^n(s)>\nonumber\\
&:=&A_{21,1}+A_{21,2}+A_{21,3}.
\end{eqnarray}
By Ito isometry and H\"o{}lder's inequality,
\begin{eqnarray}\label{3.9}
 &&E[A_{21,1}]\nonumber\\
 &\leq &
C\int_0^t(E[\int_{s_n^-}^s|f_n(u)\sigma^*\nabla\phi(X(u))-f_n(s_n^-)\sigma^*\nabla\phi(X(s_n^-))|^2du])^{\frac{1}{2}}
(E[|\dot{W}^n|^2(s)])^{\frac{1}{2}}ds \nonumber\\
&\leq& C\int_0^t(2^n)^{\frac{1}{2}}(E[\int_{s_n^-}^s|f_n(u)\sigma^*\nabla\phi(X(u))-f_n(s_n^-)\sigma^*\nabla\phi(X(s_n^-))|^2du])^{\frac{1}{2}}
ds \nonumber\\
&\leq& C\int_0^t(2^n)^{\frac{1}{2}}(\frac{1}{2^n})^{\frac{1}{2}}(\frac{1}{2^n})^{\frac{1}{2}}ds\leq C(\frac{1}{2^n})^{\frac{1}{2}},
\end{eqnarray}
where (\ref{3.01}), (\ref{3.02}) have been used. For the term $A_{21,2}$, we have
\begin{eqnarray}\label{3.10}
 &&E[A_{21,2}]\nonumber\\
 &\leq &
C\int_0^t(E[|W(s)-W(s_n^-)|^3])^{\frac{1}{3}}
(E[|\sigma^*\nabla\phi(X^n(s))-\sigma^*\nabla\phi(X^n(s_n^-))|^3])^{\frac{1}{3}}\nonumber\\
&&\quad\quad \times
(E[|\dot{W}^n|^3(s)])^{\frac{1}{3}}ds \nonumber\\
&\leq& C\int_0^t(2^n)^{\frac{1}{2}}(\frac{1}{2^n})^{\frac{1}{2}}(\frac{1}{2^n})^{\frac{1}{2}}ds\leq C(\frac{1}{2^n})^{\frac{1}{2}}.
\end{eqnarray}
where (\ref{3.01}) has been used.
Now,
\begin{eqnarray}\label{3.11}
 &&A_{21,3}\nonumber\\
 &=&r^2\sum_{k} \int_{\frac{k}{2^n}}^{\frac{k+1}{2^n}}f_n({\frac{k-1}{2^n}})<\nabla\phi(X({\frac{k-1}{2^n}})), \sigma(X({\frac{k-1}{2^n}})(W(s)-W({\frac{k}{2^n}}))>\nonumber\\
 &&\times <\nabla\phi(X^n({\frac{k-1}{2^n}})), \sigma(X^n({\frac{k-1}{2^n}}))(W({\frac{k}{2^n}})-W({\frac{k-1}{2^n}}))>ds\nonumber\\
 &+&r^2\sum_{k} \int_{\frac{k}{2^n}}^{\frac{k+1}{2^n}}f_n({\frac{k-1}{2^n}})<\sigma^*\nabla\phi(X({\frac{k-1}{2^n}})), W({\frac{k}{2^n}})-W({\frac{k-1}{2^n}}))>\nonumber\\
 &&\times <\sigma^*\nabla\phi(X^n({\frac{k-1}{2^n}})), W({\frac{k}{2^n}})-W({\frac{k-1}{2^n}})>ds\nonumber\\
 &:=&A_{21,31}+A_{21,32}.
\end{eqnarray}
Conditioning on ${\cal F}_{\frac{k}{2^n}}$, it is easy to see that $E[A_{21,31}]=0$. Moreover,
 \begin{eqnarray}\label{3.12}
 &&A_{21,32}\nonumber\\
 &=&r^2\sum_{k} f_n({\frac{k-1}{2^n}})\sum_{i=1}^m (\sigma^*\nabla\phi)_i(X({\frac{k-1}{2^n}})) (\sigma^*\nabla\phi)_i(X^n({\frac{k-1}{2^n}}))\nonumber\\
 &&\times (|W_i({\frac{k}{2^n}})-W_i({\frac{k-1}{2^n}})|^2-\frac{1}{2^n})\nonumber\\
 &+&r^2\sum_{k} f_n({\frac{k-1}{2^n}})\sum_{i\not= j}(\sigma^*\nabla\phi)_i(X({\frac{k-1}{2^n}})) (\sigma^*\nabla\phi)_j(X^n({\frac{k-1}{2^n}}))\nonumber\\
 &&\times (W_i({\frac{k}{2^n}})-W_i({\frac{k-1}{2^n}}))(W_j({\frac{k}{2^n}})-W_j({\frac{k-1}{2^n}}))\nonumber\\
 &+&r^2\sum_{k} f_n({\frac{k-1}{2^n}})\sum_{i=1}^m (\sigma^*\nabla\phi)_i(X({\frac{k-1}{2^n}})) (\sigma^*\nabla\phi)_i(X^n({\frac{k-1}{2^n}}))(\frac{1}{2^n})\nonumber\\
 &:=& A_{21,321}+A_{21,322}+A_{21,323}.
\end{eqnarray}
Conditioning on ${\cal F}_{\frac{k-1}{2^n}}$ and using the independence of $W_i, W_j$ for $i\not = j$, we find that $E[A_{21,321}]=0$ and  $E[A_{21,322}]=0$.
On the other hand,
\begin{eqnarray}\label{3.13}
 &&E[A_{21,323}]\nonumber\\
 &=&r^2E[\int_0^tf_n(s)<\sigma^*\nabla\phi(X(s)),  \sigma^*\nabla\phi(X^n(s))>ds]\nonumber\\
 &+&r^2E[\int_0^t\{f_n(s_n^-)<\sigma^*\nabla\phi(X(s_n^-)),  \sigma^*\nabla\phi(X^n(s_n^-))>\nonumber\\
 && -f_n(s)<\sigma^*\nabla\phi(X(s)),  \sigma^*\nabla\phi(X^n(s))>\} ds]\nonumber\\
 &\leq& r^2E[\int_0^tf_n(s)<\sigma^*\nabla\phi(X(s)),  \sigma^*\nabla\phi(X^n(s))> ds]+C(\frac{1}{2^n})^{\frac{1}{2}},
\end{eqnarray}
where (\ref{3.01}), (\ref{3.02}) again have been used.
Putting together (\ref{3.8})---(\ref{3.13}) we arrive at
\begin{equation}\label{3.14}
 E[A_{21}]\leq CE[\int_0^tf_n(s)ds]+C(\frac{1}{2^n})^{\frac{1}{2}}.
\end{equation}
The term $A_{25}$ can be bounded as follows.
\begin{eqnarray}\label{3.15}
 &&E[A_{25}]\nonumber\\
 &\leq &C E[\sum_{k} \int_{\frac{k}{2^n}}^{\frac{k+1}{2^n}}(\int_{\frac{k-1}{2^n}}^{s}d|L|(u))2^n|W(\frac{k}{2^n})-W(\frac{k-1}{2^n})|ds]\nonumber\\
 &\leq& C E[\sum_{k} (|L|({\frac{k}{2^n}})-|L|({\frac{k-1}{2^n}}))|W(\frac{k}{2^n})-W(\frac{k-1}{2^n})|]\nonumber\\
 &\leq& 2CE[|L|(t) \sup_{|u-v|\leq \frac{1}{2^n}}(|W(u)-W(v)|)]\nonumber\\
 &\leq& 2C(E[|L|^2(t)])^{\frac{1}{2}}(\frac{1}{2^n})^{\frac{1}{2}}\leq C(\frac{1}{2^n})^{\frac{1}{2}}.
\end{eqnarray}
To control the term $A_{26}$, we write it as
\begin{eqnarray}\label{3.16}
 &&A_{26}\nonumber\\
 &=& r^2\sum_k \int_{\frac{k}{2^n}}^{\frac{k+1}{2^n}}(\int_{\frac{k-1}{2^n}}^{\frac{k}{2^n}}f_n(u)<\nabla\phi(X^n(u)), \sigma(X^n(u))dW^n(u)> )\nonumber\\
 &&\times <\nabla\phi(X^n(s)), \sigma(X^n(s))dW^n(s)>\nonumber\\
 &+&r^2\sum_k \int_{\frac{k}{2^n}}^{\frac{k+1}{2^n}}(\int_{\frac{k}{2^n}}^{s}f_n(u)<\nabla\phi(X^n(u)), \sigma(X^n(u))dW^n(u)> )\nonumber\\
 &&\times <\nabla\phi(X^n(s)), \sigma(X^n(s))dW^n(s)>\nonumber\\
 &=&r^2\sum_k \int_{\frac{k}{2^n}}^{\frac{k+1}{2^n}}(\int_{\frac{k-1}{2^n}}^{\frac{k}{2^n}}<f_n(u)\sigma^*\nabla\phi(X^n(u))-f_n(\frac{k-1}{2^n})\sigma^*\nabla\phi(X^n(\frac{k-1}{2^n})), \nonumber\\
 && dW^n(u)> )<\nabla\phi(X^n(s)), \sigma(X^n(s))dW^n(s)>\nonumber\\
&+& r^2\sum_k \int_{\frac{k}{2^n}}^{\frac{k+1}{2^n}}(\int_{\frac{k-1}{2^n}}^{\frac{k}{2^n}}f_n(\frac{k-1}{2^n})<\sigma^*\nabla\phi(X^n(\frac{k-1}{2^n})), dW^n(u)>) \nonumber\\
&&\times <\sigma^*\nabla\phi(X^n(s))-\sigma^*\nabla\phi(X^n(\frac{k-1}{2^n})), dW^n(s)>\nonumber\\
&+& r^2\sum_k \int_{\frac{k}{2^n}}^{\frac{k+1}{2^n}}(\int_{\frac{k-1}{2^n}}^{\frac{k}{2^n}}f_n(\frac{k-1}{2^n})<\sigma^*\nabla\phi(X^n(\frac{k-1}{2^n})), dW^n(u)> )\nonumber\\
&&\times <\sigma^*\nabla\phi(X^n(\frac{k-1}{2^n})), dW^n(s)>
\nonumber\\
 &+&r^2\sum_k \int_{\frac{k}{2^n}}^{\frac{k+1}{2^n}}(\int_{\frac{k}{2^n}}^{s}f_n(u)<\nabla\phi(X^n(u)), \sigma(X^n(u))dW^n(u)> )\nonumber\\
 &&\times <\nabla\phi(X^n(s)), \sigma(X^n(s))dW^n(s)>\nonumber\\
 &:=& A_{26,1}+A_{26,2}+A_{26,3}+A_{26,4}
\end{eqnarray}
The first term on the right can be bounded as follows:
\begin{eqnarray}\label{3.17}
 &&E[A_{26,1}]\nonumber\\
 &\leq &C\sum_k \int_{\frac{k}{2^n}}^{\frac{k+1}{2^n}}ds\int_{\frac{k-1}{2^n}}^{\frac{k}{2^n}}du E[|\dot{W}^n(u)| |\dot{W}^n(s)|\nonumber\\
 &&\times |f_n(u)\sigma^*\nabla\phi(X^n(u))-f_n(\frac{k-1}{2^n})\sigma^*\nabla\phi(X^n(\frac{k-1}{2^n}))|]\nonumber\\
 &\leq &C\sum_k \int_{\frac{k}{2^n}}^{\frac{k+1}{2^n}}ds\int_{\frac{k-1}{2^n}}^{\frac{k}{2^n}}du (2^n)^{\frac{1}{2}}(2^n)^{\frac{1}{2}} \nonumber\\
 &&\times (E[|f_n(u)\sigma^*\nabla\phi(X^n(u))-f_n(\frac{k-1}{2^n})\sigma^*\nabla\phi(X^n(\frac{k-1}{2^n}))|^3])^{\frac{1}{3}}\nonumber\\
&\leq & C\sum_k \int_{\frac{k}{2^n}}^{\frac{k+1}{2^n}}ds\int_{\frac{k-1}{2^n}}^{\frac{k}{2^n}}du (2^n)^{\frac{1}{2}}(2^n)^{\frac{1}{2}} (\frac{1}{2^n})^{\frac{1}{2}}\leq C (\frac{1}{2^n})^{\frac{1}{2}}.
\end{eqnarray}
The second term has the following upper bound.
\begin{eqnarray}\label{3.18}
 &&E[A_{26,2}]\nonumber\\
 &\leq& C\sum_k \int_{\frac{k}{2^n}}^{\frac{k+1}{2^n}}ds\int_{\frac{k-1}{2^n}\vee 0}^{\frac{k}{2^n}}du E[|\dot{W}^n(u)| |\dot{W}^n(s)||\sigma^*\nabla\phi(X^n(s))-\sigma^*\nabla\phi(X^n(\frac{k-1}{2^n}))|]\nonumber\\
 &\leq &C\sum_k \int_{\frac{k}{2^n}}^{\frac{k+1}{2^n}}ds\int_{\frac{k-1}{2^n}\vee 0}^{\frac{k}{2^n}}du (2^n)^{\frac{1}{2}}(2^n)^{\frac{1}{2}} (E[|\sigma^*\nabla\phi(X^n(s))-\sigma^*\nabla\phi(X^n(\frac{k-1}{2^n}))|^3])^{\frac{1}{3}}\nonumber\\
&\leq & C\sum_k \int_{\frac{k}{2^n}}^{\frac{k+1}{2^n}}ds\int_{\frac{k-1}{2^n}\vee 0}^{\frac{k}{2^n}}du (2^n)^{\frac{1}{2}}(2^n)^{\frac{1}{2}} (\frac{1}{2^n})^{\frac{1}{2}}\leq C (\frac{1}{2^n})^{\frac{1}{2}}.
\end{eqnarray}
Note that
\begin{eqnarray}\label{3.19}
 &&A_{26,3}\nonumber\\
 &=&r^2\sum_k (2^n)^2\int_{\frac{k}{2^n}}^{\frac{k+1}{2^n}}ds\int_{\frac{k-1}{2^n}}^{\frac{k}{2^n}}du f_n(\frac{k-1}{2^n})\sum_{i,j=1}^m (\sigma^*\nabla\phi)_i(X^n(\frac{k-1}{2^n})(\sigma^*\nabla\phi)_j(X^n(\frac{k-1}{2^n})))\nonumber\\
 &&\times (W_i(\frac{k-1}{2^n})-W_i(\frac{k-2}{2^n}))(W_j(\frac{k}{2^n})-W_j(\frac{k-1}{2^n})).
\end{eqnarray}
Conditioning on ${\cal F}_{\frac{k-1}{2^n}}$ we see that $E[A_{26,3}]=0$. On the other hand,  the term  $A_{26,4}$ can be further split as follows.
\begin{eqnarray}\label{3.20}
 &&A_{26,4}\nonumber\\
 &=& r^2\sum_k \int_{\frac{k}{2^n}}^{\frac{k+1}{2^n}}(\int_{\frac{k}{2^n}}^{s}
 <f_n(u)\sigma^*\nabla\phi(X^n(u))-f_n(\frac{k-1}{2^n})\sigma^*\nabla\phi(X^n(\frac{k-1}{2^n})), \dot{W}^n(u)>du )\nonumber\\
 &&\times <\sigma^*\nabla\phi(X^n(s)), \dot{W}^n(s)>ds\nonumber\\
 &+& r^2\sum_k \int_{\frac{k}{2^n}}^{\frac{k+1}{2^n}}(\int_{\frac{k}{2^n}}^{s}
 <f_n(\frac{k-1}{2^n})\sigma^*\nabla\phi(X^n(\frac{k-1}{2^n})), \dot{W}^n(u)>du )\nonumber\\
 &&\times <\sigma^*\nabla\phi(X^n(s))-\sigma^*\nabla\phi(X^n(\frac{k-1}{2^n})), \dot{W}^n(s)>ds\nonumber\\
 &+& r^2(2^n)^2\sum_k \int_{\frac{k}{2^n}}^{\frac{k+1}{2^n}}ds\int_{\frac{k}{2^n}}^{s}du
 f_n(\frac{k-1}{2^n})<\sigma^*\nabla\phi(X^n(\frac{k-1}{2^n})), W(\frac{k}{2^n})-W(\frac{k-1}{2^n})>^2\nonumber\\
 &:=&A_{26,41}+A_{26,42}+A_{26,43}
\end{eqnarray}
By the Lipschitz continuity of the coefficients and (\ref{3.01}) and (\ref{3.02})  we have
\begin{eqnarray}\label{3.21}
 &&E[A_{26,41}]\nonumber\\
 &\leq &C\sum_k \int_{\frac{k}{2^n}}^{\frac{k+1}{2^n}}ds\int_{\frac{k}{2^n}}^{s}du E[|\dot{W}^n(u)| |\dot{W}^n(s)|\nonumber\\
  &&\times |f_n(u)\sigma^*\nabla\phi(X^n(u))-f_n(\frac{k-1}{2^n})\sigma^*\nabla\phi(X^n(\frac{k-1}{2^n}))|]\nonumber\\
&\leq & C\sum_k \int_{\frac{k}{2^n}}^{\frac{k+1}{2^n}}ds\int_{\frac{k}{2^n}}^{\frac{k+1}{2^n}}du (2^n)^{\frac{1}{2}}(2^n)^{\frac{1}{2}} (\frac{1}{2^n})^{\frac{1}{2}}\leq C (\frac{1}{2^n})^{\frac{1}{2}}.
\end{eqnarray}
and
\begin{eqnarray}\label{3.22}
 &&E[A_{26,42}]\nonumber\\
 &\leq&C\sum_k \int_{\frac{k}{2^n}}^{\frac{k+1}{2^n}}ds\int_{\frac{k}{2^n}}^{s}du E[|\dot{W}^n(u)| |\dot{W}^n(s)| \cdot |\sigma^*\nabla\phi(X^n(s))-\sigma^*\nabla\phi(X^n(\frac{k-1}{2^n}))|]\nonumber\\
&\leq & C\sum_k \int_{\frac{k}{2^n}}^{\frac{k+1}{2^n}}ds\int_{\frac{k}{2^n}}^{\frac{k+1}{2^n}}du (2^n)^{\frac{1}{2}}(2^n)^{\frac{1}{2}} (\frac{1}{2^n})^{\frac{1}{2}}\leq C (\frac{1}{2^n})^{\frac{1}{2}}.
\end{eqnarray}
Furthermore,
\begin{eqnarray}\label{3.23}
 &&A_{26,43}\nonumber\\
 &=& \frac{1}{2} r^2\sum_k f_n(\frac{k-1}{2^n})\sum_{i\not=j}^m (\sigma^*\nabla\phi)_i(X^n(\frac{k-1}{2^n}))(\sigma^*\nabla\phi)_j(X^n(\frac{k-1}{2^n}))\nonumber\\
 &&\times ( W_i(\frac{k}{2^n})-W_i(\frac{k-1}{2^n}))( W_j(\frac{k}{2^n})-W_j(\frac{k-1}{2^n}))\nonumber\\
 &+&\frac{1}{2} r^2\sum_k f_n(\frac{k-1}{2^n})\sum_{i=1}^m (\sigma^*\nabla\phi)_i^2(X^n(\frac{k-1}{2^n}))( |W_i(\frac{k}{2^n})-W_i(\frac{k-1}{2^n})|^2-\frac{1}{2^n})\nonumber\\
 &+&\frac{1}{2} r^2\int_0^t[f_n(s_n^-)|\sigma^*\nabla\phi|^2(X^n(s_n^-))-f_n(s)|\sigma^*\nabla\phi|^2(X^n(s))]ds \nonumber\\
 &+&\frac{1}{2} r^2\int_0^tf_n(s)|\sigma^*\nabla\phi|^2(X^n(s))ds.
\end{eqnarray}
Conditioning on ${\cal F}_{\frac{k-1}{2^n}}$ and using the independence of $W_i, W_j$ for $i\not =j$, it is easy to see that the expectation of the first two terms on the right side are zero. By (\ref{3.01}), the expectation of the third term is bounded by $C(\frac{1}{2^n})^{\frac{1}{2}}$. Thus we conclude that
\begin{eqnarray}\label{3.24}
 &&E[A_{26,43}]\nonumber\\
 &\leq & C(\frac{1}{2^n})^{\frac{1}{2}}+\frac{1}{2} r^2E[\int_0^tf_n(s)|\sigma^*\nabla\phi|^2(X^n(s))ds ].
\end{eqnarray}
Combining (\ref{3.20})---(\ref{3.24}), we find that
\begin{eqnarray}\label{3.25}
 &&E[A_{26,4}]\nonumber\\
 &\leq & C(\frac{1}{2^n})^{\frac{1}{2}}+\frac{1}{2} r^2E[\int_0^tf_n(s)|\sigma^*\nabla\phi|^2(X^n(s))ds ].
\end{eqnarray}
Putting together (\ref{3.16}),(\ref{3.17}),(\ref{3.18}), (\ref{3.19}) and (\ref{3.25}) yields
\begin{eqnarray}\label{3.26}
 &&E[A_{26}]\nonumber\\
 &\leq & C(\frac{1}{2^n})^{\frac{1}{2}}+\frac{1}{2} r^2E[\int_0^tf_n(s)|\sigma^*\nabla\phi|^2(X^n(s))ds ].
\end{eqnarray}
The term $A_{28}$ admits a similar bound as  $A_{25}$:
\begin{eqnarray}\label{3.27}
 &&E[A_{28}]\nonumber\\
 &\leq& CE[|L^n|(t) \sup_{|u-v|\leq \frac{1}{2^n}}(|W(u)-W(v)|)]\nonumber\\
 &\leq& 2C[\sup_n(E[|L^n|^2(t)])^{\frac{1}{2}}](\frac{1}{2^n})^{\frac{1}{2}}\leq C(\frac{1}{2^n})^{\frac{1}{2}}.
\end{eqnarray}
Now let us turn  to $A_{29}$. We have
\begin{eqnarray}\label{3.28}
 &&A_{29}\nonumber\\
 &=&2r\sum_{k}\int_{\frac{k}{2^n}}^{\frac{k+1}{2^n}}(\int_{\frac{k-1}{2^n}}^{\frac{k}{2^n}}g_n(u)<X^n(u)-X(u), \sigma(X^n(u))dW^n(u)>)\nonumber\\
 &&\times <\nabla\phi(X^n(s)), \sigma(X^n(s))dW^n(s)> \nonumber\\
 &+&2r\sum_{k}\int_{\frac{k}{2^n}}^{\frac{k+1}{2^n}}(\int_{\frac{k}{2^n}}^{s}g_n(u)<X^n(u)-X(u), \sigma(X^n(u))dW^n(u)>)\nonumber\\
 &&\times <\nabla\phi(X^n(s)), \sigma(X^n(s))dW^n(s)> \nonumber\\
 &=&2r\sum_{k}\int_{\frac{k}{2^n}}^{\frac{k+1}{2^n}}(\int_{\frac{k-1}{2^n}}^{\frac{k}{2^n}}g_n(u)<X^n(u)-X(u), \sigma(X^n(u))dW^n(u)>)\nonumber\\
 &&\times <\sigma^*\nabla\phi(X^n(s))-\sigma^*\nabla\phi(X^n(\frac{k-1}{2^n})), dW^n(s)> \nonumber\\
 &+&2r\sum_{k}\int_{\frac{k}{2^n}}^{\frac{k+1}{2^n}}(\int_{\frac{k-1}{2^n}}^{\frac{k}{2^n}}<g_n(u)\sigma^*(X^n(u))(X^n(u)-X(u))
 \nonumber\\
 &&-g_n(\frac{k-1}{2^n})\sigma^*(X^n(\frac{k-1}{2^n}))(X^n(\frac{k-1}{2^n})
 -X(\frac{k-1}{2^n})), dW^n(u)>)\nonumber\\
 &&\quad\quad \times <\sigma^*\nabla\phi(X^n(\frac{k-1}{2^n})), dW^n(s)> \nonumber\\
 &+&2r\sum_{k}<g_n(\frac{k-1}{2^n})\sigma^*(X^n(\frac{k-1}{2^n}))(X^n(\frac{k-1}{2^n})
 -X(\frac{k-1}{2^n})), W(\frac{k-1}{2^n})-W(\frac{k-2}{2^n})>\nonumber\\
 &&\times <\sigma^*\nabla\phi(X^n(\frac{k-1}{2^n})), W(\frac{k}{2^n})-W(\frac{k-1}{2^n})> \nonumber\\
 &+&2r\sum_{k}\int_{\frac{k}{2^n}}^{\frac{k+1}{2^n}}(\int_{\frac{k}{2^n}}^{s}g_n(u)<X^n(u)-X(u), \sigma(X^n(u))dW^n(u)>)\nonumber\\
 &&\times <\sigma^*\nabla\phi(X^n(s))-\sigma^*\nabla\phi(X^n(\frac{k-1}{2^n})), dW^n(s)> \nonumber\\
 &+&2r\sum_{k}\int_{\frac{k}{2^n}}^{\frac{k+1}{2^n}}(\int_{\frac{k}{2^n}}^{s}<g_n(u)\sigma^*(X^n(u))(X^n(u)-X(u))
 -g_n(\frac{k-1}{2^n})\sigma^*(X^n(\frac{k-1}{2^n}))\nonumber\\
 &&\times (X^n(\frac{k-1}{2^n})
 -X(\frac{k-1}{2^n})), dW^n(u)>)<\sigma^*\nabla\phi(X^n(\frac{k-1}{2^n})), dW^n(s)> \nonumber\\
 &+&r\sum_{k}<g_n(\frac{k-1}{2^n})\sigma^*(X^n(\frac{k-1}{2^n}))(X^n(\frac{k-1}{2^n})
 -X(\frac{k-1}{2^n})), W(\frac{k}{2^n})-W(\frac{k-1}{2^n})>\nonumber\\
 &&\times <\sigma^*\nabla\phi(X^n(\frac{k-1}{2^n})), W(\frac{k}{2^n})-W(\frac{k-1}{2^n})>\nonumber\\
 &:=& \sum_{i=1}^6 A_{29,i}
\end{eqnarray}
The first and the second term on the right have the following bounds.
\begin{eqnarray}\label{3.29}
 &&E[A_{29,1}]\nonumber\\
 &\leq &2r\sum_{k}\int_{\frac{k}{2^n}}^{\frac{k+1}{2^n}}ds\int_{\frac{k-1}{2^n}}^{\frac{k}{2^n}}du E[|\dot{W}^n(u)||\dot{W}^n(s)||\sigma^*\nabla\phi(X^n(s))-\sigma^*\nabla\phi(X^n(\frac{k-1}{2^n})|]\nonumber\\
 &\leq& C(\frac{1}{2^n})^{\frac{1}{2}}.
 \end{eqnarray}
 \begin{eqnarray}\label{3.30}
 &&E[A_{29,2}]\nonumber\\
 &\leq &2r\sum_{k}\int_{\frac{k}{2^n}}^{\frac{k+1}{2^n}}ds\int_{\frac{k-1}{2^n}}^{\frac{k}{2^n}}du E[|\dot{W}^n(u)||\dot{W}^n(s)|\nonumber\\
 &&\times |g_n(u)\sigma^*(X^n(u))(X^n(u)-X(u))-g_n(\frac{k-1}{2^n})\sigma^*(X^n(\frac{k-1}{2^n}))(X^n(\frac{k-1}{2^n})
 -X(\frac{k-1}{2^n}))|]\nonumber\\
 &\leq& C(\frac{1}{2^n})^{\frac{1}{2}},
 \end{eqnarray}
 here (\ref{3.01}), (\ref{3.02}) have been used again.
 By a similar reason, it also holds that
 \begin{equation}\label{3.31}
 E[A_{29,4}]\leq C(\frac{1}{2^n})^{\frac{1}{2}}.
 \end{equation}
 \begin{equation}\label{3.32}
 E[A_{29,5}]\leq C(\frac{1}{2^n})^{\frac{1}{2}}.
 \end{equation}
 Conditioning on ${\cal F}_{\frac{k-1}{2^n}}$, we have that $E[A_{29,3}]=0$.
 Note that
 \begin{eqnarray}\label{3.33}
 &&A_{29,6}\nonumber\\
 &=&r\sum_{k}\sum_{i\not=j}[g_n(\frac{k-1}{2^n})(\sigma^*(X^n(\frac{k-1}{2^n}))(X^n(\frac{k-1}{2^n})
 -X(\frac{k-1}{2^n})))_i (\sigma^*\nabla\phi)_j(X^n(\frac{k-1}{2^n}))\nonumber\\
 &&\times (W_i(\frac{k}{2^n})-W_i(\frac{k-1}{2^n}))( W_i(\frac{k}{2^n})-W_i(\frac{k-1}{2^n}))]\nonumber\\
 &+& r\sum_{k}\sum_{i=1}^m[g_n(\frac{k-1}{2^n})(\sigma^*(X^n(\frac{k-1}{2^n}))(X^n(\frac{k-1}{2^n})
 -X(\frac{k-1}{2^n})))_i (\sigma^*\nabla\phi)_i(X^n(\frac{k-1}{2^n}))\nonumber\\
 &&\times \{|W_i(\frac{k}{2^n})-W_i(\frac{k-1}{2^n})|^2-\frac{1}{2^n}\}\nonumber\\
 &+& r\int_0^t \{<g_n(s_n^-)\sigma^*(X^n(s_n^-))(X^n(s_n^-)
 -X(s_n^-)), \sigma^*\nabla\phi(X^n(s_n^-))>\nonumber\\
 && -<g_n(s)\sigma^*(X^n(s))(X^n(s)
 -X(s)), \sigma^*\nabla\phi(X^n(s))>\}ds\nonumber\\
 &+&r\int_0^t <g_n(s)\sigma^*(X^n(s))(X^n(s)
 -X(s)), \sigma^*\nabla\phi(X^n(s))>ds\nonumber\\
 &:=&A_{29,61}+A_{29,62}+A_{29,63}+A_{29,64}
\end{eqnarray}
Again by conditioning on ${\cal F}_{\frac{k-1}{2^n}}$ and the independence,
$$E[A_{29,61}]=0,\quad E[A_{29,62}]=0.$$
By virtue of (\ref{3.01}) and (\ref{3.02}),
 \begin{equation}\label{3.34}
 E[A_{29,63}]\leq C(\frac{1}{2^n})^{\frac{1}{2}}.
 \end{equation}
 It follows from (\ref{3.28})---(\ref{3.34}) that
 \begin{eqnarray}\label{3.35}
 &&E[A_{29}]\nonumber\\
 &\leq & C(\frac{1}{2^n})^{\frac{1}{2}}+ r\int_0^t <g_n(s)\sigma^*(X^n(s))(X^n(s)
 -X(s)), \sigma^*\nabla\phi(X^n(s))>ds.\nonumber\\
 &&
\end{eqnarray}
Applying the same arguments to $A_{210}$ in (\ref{3.5}), we get
\begin{eqnarray}\label{3.36}
 &&E[A_{210}]\nonumber\\
 &\leq & C(\frac{1}{2^n})^{\frac{1}{2}}- 2r\int_0^t <g_n(s)\sigma^*(X(s))(X^n(s)
 -X(s)), \sigma^*\nabla\phi(X^n(s))>ds.\nonumber\\
 &&
\end{eqnarray}
As for the term $A_{213}$ in (\ref{3.5}), we have
\begin{eqnarray}\label{3.37}
 &&E[A_{213}]\nonumber\\
 &\leq &C E[\sum_{k} \int_{\frac{k}{2^n}}^{\frac{k+1}{2^n}}(\int_{\frac{k-1}{2^n}}^{s}d|L^n|_u)(2^n|W(\frac{k}{2^n})-W(\frac{k-1}{2^n})|ds\nonumber\\
 &\leq& C E[\sum_{k} (|L^n|_{\frac{k}{2^n}}-|L^n|_{\frac{k-1}{2^n}})|W(\frac{k}{2^n})-W(\frac{k-1}{2^n})|\nonumber\\
 &\leq& 2CE[|L^n|_t \sup_{|u-v|\leq \frac{1}{2^n}}(|W(u)-W(v)|)]\nonumber\\
 &\leq& 2C(E[|L^n|_t^2])^{\frac{1}{2}}(\frac{1}{2^n})^{\frac{1}{2}}\leq C(\frac{1}{2^n})^{\frac{1}{2}}.
\end{eqnarray}
A similar argument leads to
\begin{eqnarray}\label{3.38}
 E[A_{214}] &\leq&  C(\frac{1}{2^n})^{\frac{1}{2}}.
\end{eqnarray}
Collecting the estimates (\ref{3.5})---(\ref{3.38}) we  get that
\begin{eqnarray}\label{3.39}
 &&E[A_2] \nonumber\\
 &\leq&  C(\frac{1}{2^n})^{\frac{1}{2}}+r^2E[\int_0^tf_n(s)<\sigma^*\nabla\phi(X^n(s)),\sigma^*\nabla\phi(X(s))> ds]\nonumber\\
 &&+\frac{1}{2} r^2E[\int_0^tf_n(s)|\sigma^*\nabla\phi|^2(X^n(s))ds\nonumber\\
 &&+ r\int_0^t <g_n(s)\sigma^*(X^n(s))(X^n(s)
 -X(s)), \sigma^*\nabla\phi(X^n(s))>ds\nonumber\\
 &&-2r\int_0^t <g_n(s)\sigma^*(X(s))(X^n(s)
 -X(s)), \sigma^*\nabla\phi(X^n(s))>ds.\nonumber\\
 &&
\end{eqnarray}
Now we turn to $A_3$. By the chain rule, we have
\begin{eqnarray}\label{3.40}
 &&A_3\nonumber\\
  &=& r\int_0^tf_n(s_n^{-})\sum_{i=1}^m[(\sigma^*\nabla\phi)_i(X^n(s))-(\sigma^*\nabla\phi)_i(X^n(s_n^{-}))]dW_i^n(s)\nonumber\\
&=& r\int_0^tf_n(s_n^{-})\sum_{i=1}^m\int_{s_n^-}^s[<\nabla(\sigma^*\nabla\phi)_i(X^n(u))
-\nabla(\sigma^*\nabla\phi)_i(X^n(s_n^-)), \nonumber\\
&&\quad\quad \sigma(X^n(u))dW^n(u)>]dW_i^n(s)\nonumber\\
&+&r\int_0^tf_n(s_n^{-})\sum_{i=1}^m\int_{s_n^-}^s<\nabla(\sigma^*\nabla\phi)_i(X^n(s_n^-)), (\sigma(X^n(u))\nonumber\\
&&\quad\quad -\sigma(X^n(s_n^-)))dW^n(u)>dW_i^n(s)\nonumber\\
&+&r\int_0^tf_n(s_n^{-})\sum_{i=1}^m\int_{s_n^-}^s<\nabla(\sigma^*\nabla\phi)_i(X^n(s_n^-)), \sigma(X^n(s_n^-))dW^n(u)>dW_i^n(s)\nonumber\\
&+&r\int_0^tf_n(s_n^{-})\sum_{i=1}^m\int_{s_n^-}^s<\nabla(\sigma^*\nabla\phi)_i(X^n(u)), \nu(X^n(u))d|L^n|(u)>dW_i^n(s)\nonumber\\
&+& r\int_0^tf_n(s_n^{-})\sum_{i=1}^m\int_{s_n^-}^s<\nabla(\sigma^*\nabla\phi)_i(X^n(u)), b(X^n(u))du>dW_i^n(s)\nonumber\\
&:=& A_{31}+A_{32}+A_{33}+A_{34}+A_{35}
\end{eqnarray}
Similar to the estimates for $A_{214}$, $A_{22}$ and the term $A_{21,2}$, it can be shown that
\begin{equation}\label{3.41}
E[A_{3i}]\leq  C(\frac{1}{2^n})^{\frac{1}{2}}, \quad i=1,2,4,5.
\end{equation}
Now,
\begin{eqnarray}\label{3.42}
 &&A_{33}\nonumber\\
  &=& r\sum_{k}(2^n)^2\int_{\frac{k}{2^n}}^{\frac{k+1}{2^n}}ds\int_{\frac{k-1}{2^n}}^{\frac{k}{2^n}}du f_n(\frac{k-1}{2^n})\sum_{i=1}^m\sum_{j=1}^m(\sigma^*(\nabla(\sigma^*\nabla\phi)_i))_j(X^n(\frac{k-1}{2^n}))
  \nonumber\\
  &&\times (W_i(\frac{k}{2^n})-W_i(\frac{k-1}{2^n}))(W_j(\frac{k-1}{2^n}-W_j(\frac{k-2}{2^n}))\nonumber\\
&+&r\sum_{k}(2^n)^2\int_{\frac{k}{2^n}}^{\frac{k+1}{2^n}}ds\int_{\frac{k}{2^n}}^{s}du f_n(\frac{k-1}{2^n})\sum_{i=1}^m\sum_{j=1}^m(\sigma^*(\nabla(\sigma^*\nabla\phi)_i))_j(X^n(\frac{k-1}{2^n}))
  \nonumber\\
  &&\times (W_i(\frac{k}{2^n})-W_i(\frac{k-1}{2^n}))(W_j(\frac{k}{2^n})-W_j(\frac{k-1}{2^n}))\nonumber\\
  &=&r\sum_{k}f_n(\frac{k-1}{2^n})\sum_{i=1}^m\sum_{j=1}^m(\sigma^*(\nabla(\sigma^*\nabla\phi)_i))_j(X^n(\frac{k-1}{2^n}))
  \nonumber\\
  &&\times (W_i(\frac{k}{2^n})-W_i(\frac{k-1}{2^n}))(W_j(\frac{k-1}{2^n})-W_j(\frac{k-2}{2^n}))\nonumber\\
&+&\frac{1}{2}r\sum_{k}f_n(\frac{k-1}{2^n})\sum_{i=1}^m\sum_{j=1}^m(\sigma^*(\nabla(\sigma^*\nabla\phi)_i))_j(X^n(\frac{k-1}{2^n}))
  \nonumber\\
  &&\times (W_i(\frac{k}{2^n})-W_i(\frac{k-1}{2^n}))(W_j(\frac{k}{2^n})-W_j(\frac{k-1}{2^n}))\nonumber\\
  &:=&A_{331}+A_{332}
\end{eqnarray}
Conditioning on ${\cal F}_{\frac{k-1}{2^n}}$, it is easy to see $E[A_{331}]=0$. For the second term  we  have
\begin{eqnarray}\label{3.43}
 &&A_{332}\nonumber\\
  &=&\frac{1}{2}r\sum_{k}f_n(\frac{k-1}{2^n})\sum_{i\not=j}^m(\sigma^*(\nabla(\sigma^*\nabla\phi)_i))_j(X^n(\frac{k-1}{2^n}))
  \nonumber\\
  &&\times (W_i(\frac{k}{2^n})-W_i(\frac{k-1}{2^n}))(W_j(\frac{k}{2^n})-W_j(\frac{k-1}{2^n}))\nonumber\\
  &+&\frac{1}{2}r\sum_{k}f_n(\frac{k-1}{2^n})\sum_{i=1}^m(\sigma^*(\nabla(\sigma^*\nabla\phi)_i))_i(X^n(\frac{k-1}{2^n}))
  \{|W_i(\frac{k}{2^n})-W_i(\frac{k-1}{2^n})|^2- \frac{1}{2^n}\} \nonumber\\
  &+&\frac{1}{2}r\int_0^t\{f_n(s_n^-)\sum_{i=1}^m(\sigma^*(\nabla(\sigma^*\nabla\phi)_i))_i(X^n(s_n^-))
  -f_n(s)\sum_{i=1}^m(\sigma^*(\nabla(\sigma^*\nabla\phi)_i))_i(X^n(s)) \}ds \nonumber\\
  &+&\frac{1}{2}r\int_0^tf_n(s)\sum_{i=1}^m(\sigma^*(\nabla(\sigma^*\nabla\phi)_i))_i(X^n(s)) ds\nonumber\\
  &:=&A_{3321}+A_{3322}+A_{3323}+A_{3324}
\end{eqnarray}
Using the  martingale property and the independence of $W_i, W_j$ for $i\not = j$, we find that $E[A_{3321}]=0$ and $E[A_{3322}]=0$.  In view of (\ref{3.01}) and (\ref{3.02}),
we have $E[A_{3323}]\leq C(\frac{1}{2^n})^{\frac{1}{2}}$. Thus, we deduce from (\ref{3.42}), (\ref{3.43}) that
\begin{eqnarray}\label{3.44}
 &&E[A_{33}]\nonumber\\
  &\leq &C(\frac{1}{2^n})^{\frac{1}{2}}+\frac{1}{2}r\int_0^tf_n(s)\sum_{i=1}^m(\sigma^*(\nabla(\sigma^*\nabla\phi)_i))_i(X^n(s)) ds.
\end{eqnarray}
Finally it follows from (\ref{3.40}), (\ref{3.41}), (\ref{3.41}) that
\begin{eqnarray}\label{3.45}
 &&E[A_3]\nonumber\\
  &\leq &C(\frac{1}{2^n})^{\frac{1}{2}}+\frac{1}{2}r\int_0^tf_n(s)\sum_{i=1}^m(\sigma^*(\nabla(\sigma^*\nabla\phi)_i))_i(X^n(s)) ds.
\end{eqnarray}
Combining (\ref{3.39}) with (\ref{3.45}), we complete the proof of Lemma.
\vskip 0.3cm
\begin{lemma}
We have
\begin{eqnarray}\label{3.46}
&& rE[\int_0^tg_n(s)<X^n(s)-X(s), \sigma(X^n(s))dW^n(s)>]\nonumber\\
&\leq& rE[\int_0^tg_n(s)<\sigma^*\nabla\phi(X^n(s)),  \sigma^*(X^n(s))(X^n(s)-X(s))>ds]\nonumber\\
 &+&2rE[\int_0^tg_n(s)<\sigma^*\nabla\phi(X(s)),  \sigma^*(X^n(s))(X^n(s)-X(s))>ds]\nonumber\\
 &+&E[\int_0^tg_n(s)\sum_{i=1}^d\sum_{j=1}^m\sigma_{ij}^2(X^n(s))ds]\nonumber\\
 &+&
 E[\int_0^tg_n(s)\sum_{i=1}^d(X_i^n(s)-X_i(s))\sum_{j =1}^m (\sigma^*\nabla\sigma_{ij})_j(X^n(s))\}ds]\nonumber\\
 &-&2E[\int_0^tg_n(s)\sum_{i=1}^d\sum_{j=1}^m
 \sigma_{ij}(X(s)\sigma_{ij}(X^n(s))ds]+C(\frac{1}{2^n})^{\frac{1}{2}}.
\end{eqnarray}
\end{lemma}
\vskip 0.3cm
\noindent{\bf Proof}. Set
$$B=2\int_0^tg_n(s)<X^n(s)-X(s), \sigma(X^n(s))dW^n(s)>.$$
and write
\begin{eqnarray}\label{3.47}
B&=& 2\int_0^tg_n(s_n^{-})<X^n(s_n^{-})-X(s_n^-), \sigma(X^n(s_n^{-}))dW^n(s)>\nonumber\\
&+&2\int_0^t(g_n(s)-g_n(s_n^{-}))<X^n(s)-X(s), \sigma(X^n(s))dW^n(s)>\nonumber\\
&+&2 \int_0^tg_n(s_n^{-})<(X^n(s)-X(s))-(X^n(s_n^{-})-X(s_n^{-})),\sigma(X^n(s))dW^n(s)> \nonumber\\
&+&2 \int_0^tg_n(s_n^{-})<X^n(s_n^{-})-X(s_n^{-}),(\sigma(X^n(s))-\sigma(X^n(s_n^{-})))dW^n(s)> \nonumber\\
&:=& B_1+B_2+B_3+B_4.
\end{eqnarray}
As a stochastic integral against Brownian motion, it is easily seen that $E[B_1]=0$. In view of  (\ref{3.2}),
\begin{eqnarray}\label{3.48}
&&B_2\nonumber\\
&=& 2r\int_0^t(\int_{s_n^-}^sg_n(u)<\nabla\phi(X(u)), \sigma(X(u))dW(u)>)<X^n(s)-X(s), \sigma(X^n(s))dW^n(s)> \nonumber\\
&+&2r\int_0^t(\int_{s_n^-}^sg_n(u)<\nabla\phi(X(u)), b(X(u)du>)<X^n(s)-X(s), \sigma(X^n(s))dW^n(s)> \nonumber\\
&+&r\int_0^t(\int_{s_n^-}^sg_n(u)tr(\phi^{\prime\prime}(\sigma\sigma^*)(X(u))du)<X^n(s)-X(s), \sigma(X^n(s))dW^n(s)> \nonumber\\
&+&r\int_0^t(\int_{s_n^-}^sg_n(u)<\nabla \phi (X(u)), (\sigma\sigma^{\prime})(X(u))>du)<X^n(s)-X(s), \sigma(X^n(s))dW^n(s)> \nonumber\\
&+&2r\int_0^t(\int_{s_n^-}^sg_n(u)<\nabla \phi (X(u)), \nu(X(u))>d|L|(u))<X^n(s)-X(s), \sigma(X^n(s))dW^n(s)> \nonumber\\
&+&2r\int_0^t(\int_{s_n^-}^sg_n(u)<\nabla\phi(X^n(u)), \sigma(X^n(u))dW^n(u)> )<X^n(s)-X(s), \sigma(X^n(s))dW^n(s)> \nonumber \\
&+&2r\int_0^t(\int_{s_n^-}^sg_n(u)<\nabla\phi(X^n(u)), b(X^n(u))du>)<X^n(s)-X(s), \sigma(X^n(s))dW^n(s)> \nonumber\\
&+&2r\int_0^t(\int_{s_n^-}^sg_n(u)<\nabla \phi (X^n(u)), \nu(X^n(u))>d|L^n|(u))<X^n(s)-X(s), \sigma(X^n(s))dW^n(s)> \nonumber\\
&+&r^2\int_0^t(\int_{s_n^-}^sg_n(u)|\sigma^*\nabla\phi|^2(X(u))du)<X^n(s)-X(s), \sigma(X^n(s))dW^n(s)> \nonumber\\
&:=& \sum_{i=1}^{9}B_{2i}
\end{eqnarray}
We will closely study each of the terms on the right side. Since $\nabla\phi$, $b$, $\sigma$ are bounded on $\bar{D}$, we have
\begin{eqnarray}\label{3.49}
E[|B_{22}|]&\leq& C\int_0^t (s-s_n^-)E[|\dot{W}^n(s)|]ds\nonumber\\
&\leq& C\frac{1}{2^n}\int_0^t (2^n)^{\frac{1}{2}}ds\leq C (\frac{1}{2^n})^{\frac{1}{2}}.
\end{eqnarray}
Similar arguments lead to
\begin{eqnarray}\label{3.50}
E[|B_{2i}|]&\leq& C (\frac{1}{2^n})^{\frac{1}{2}}, \quad i=3, 4,7,9.
\end{eqnarray}
Regarding $B_{25}$, we have
\begin{eqnarray}\label{3.51}
 &&E[B_{25}]\nonumber\\
 &\leq &C E[\sum_{k} \int_{\frac{k}{2^n}}^{\frac{k+1}{2^n}}(\int_{\frac{k-1}{2^n}}^{s}d|L|(u))2^n|W(\frac{k}{2^n})-W(\frac{k-1}{2^n})|ds]\nonumber\\
 &\leq& C E[\sum_{k} (|L|({\frac{k+1}{2^n}})-|L|({\frac{k-1}{2^n}}))|W(\frac{k}{2^n})-W(\frac{k-1}{2^n})|]\nonumber\\
 &\leq& 2CE[|L|(t) \sup_{|u-v|\leq \frac{1}{2^n}}(|W(u)-W(v)|)]\nonumber\\
 &\leq& 2C(E[|L|^2(t)])^{\frac{1}{2}}(\frac{1}{2^n})^{\frac{1}{2}}\leq C(\frac{1}{2^n})^{\frac{1}{2}}.
\end{eqnarray}
Similarly,
\begin{equation}\label{3.52}
E[B_{28}]\leq C(\frac{1}{2^n})^{\frac{1}{2}}.
\end{equation}
Now,
 \begin{eqnarray}\label{3.53}
 &&B_{21}\nonumber\\
 &=&
2r\int_0^t(\int_{s_n^-}^s[g_n(u)<\nabla\phi(X(u)), \sigma(X(u)dW(u)>-g_n(s_n^-)<\nabla\phi(X(s_n^-)),\nonumber\\
&&\quad\quad  \sigma(X(s_n^-))dW(u)>])<X^n(s)-X(s)
, \sigma(X^n(s))dW^n(s)> \nonumber\\
&+& 2r\int_0^t(\int_{s_n^-}^sg_n(s_n^-)<\nabla\phi(X(s_n^-)), \sigma(X(s_n^-))dW(u)>)\nonumber\\
&&\times [<\sigma^*(X^n(s))(X^n(s)-X(s))-\sigma^*(X^n(s_n^-))(X^n(s_n^-)-X(s_n^-)), dW^n(s)>]\nonumber\\
&+& 2r\int_0^tg_n(s_n^-)<\nabla\phi(X(s_n^-)), \sigma(X(s_n^-))(W(s)-W(s_n^-))>\nonumber\\
&&\quad\quad \times <\sigma^*(X^n(s_n^-))(X^n(s_n^-)-X(s_n^-)), dW^n(s)>\nonumber\\
&:=&B_{211}+B_{212}+B_{213}.
\end{eqnarray}
By Ito isometry and H\"o{}lder's inequality,
\begin{eqnarray}\label{3.54}
 &&E[B_{211}]\nonumber\\
 &\leq &
C\int_0^t(E[\int_{s_n^-}^s|g_n(u)\sigma^*\nabla\phi(X(u))-g_n(s_n^-)\sigma^*\nabla\phi(X(s_n^-))|^2du])^{\frac{1}{2}}
(E[|\dot{W}^n|^2(s)])^{\frac{1}{2}}ds \nonumber\\
&\leq& C\int_0^t(2^n)^{\frac{1}{2}}(\frac{1}{2^n})^{\frac{1}{2}}(\frac{1}{2^n})^{\frac{1}{2}}ds\leq C(\frac{1}{2^n})^{\frac{1}{2}},
\end{eqnarray}
where (\ref{3.01}) and (\ref{3.02}) have been used. Term $B_{212}$ has the following bound.
\begin{eqnarray}\label{3.55}
 &&E[B_{212}]\nonumber\\
 &\leq &
C\int_0^t(E[|W(s)-W(s_n^-)|^3])^{\frac{1}{3}}(E[|\dot{W}^n|^3(s)])^{\frac{1}{3}}\nonumber\\
&&\times (E[|\sigma^*(X^n(s))(X^n(s)-X(s))-\sigma^*(X^n(s_n^-))(X^n(s_n^{-})-X(s_n^{-}))|^3])^{\frac{1}{3}}ds \nonumber\\
&\leq& C\int_0^t(2^n)^{\frac{1}{2}}(\frac{1}{2^n})^{\frac{1}{2}}(\frac{1}{2^n})^{\frac{1}{2}}ds\leq C(\frac{1}{2^n})^{\frac{1}{2}}.
\end{eqnarray}
Now,
\begin{eqnarray}\label{3.56}
 &&B_{213}\nonumber\\
 &=&2r\sum_{k}2^n \int_{\frac{k}{2^n}}^{\frac{k+1}{2^n}}g_n({\frac{k-1}{2^n}})<\nabla\phi(X({\frac{k-1}{2^n}})), \sigma(X({\frac{k-1}{2^n}})(W(s)-W({\frac{k}{2^n}}))>\nonumber\\
 &&\times <\sigma^*(X^n({\frac{k-1}{2^n}}))(X^n({\frac{k-1}{2^n}})-X({\frac{k-1}{2^n}})), W({\frac{k}{2^n}})-W({\frac{k-1}{2^n}})>ds\nonumber\\
 &+&2r\sum_kg_n({\frac{k-1}{2^n}})<\sigma^*\nabla\phi(X({\frac{k-1}{2^n}})), W({\frac{k}{2^n}})-W({\frac{k-1}{2^n}})>\nonumber\\
 &&\times <\sigma^*(X^n({\frac{k-1}{2^n}}))(X^n({\frac{k-1}{2^n}})-X({\frac{k-1}{2^n}})), W({\frac{k}{2^n}})-W({\frac{k-1}{2^n}})>\nonumber\\
 &:=&B_{2131}+B_{2132}.
\end{eqnarray}
Conditioning on ${\cal F}_{\frac{k}{2^n}}$, it is easy to see that $E[B_{2131}]=0$. Moreover,
 \begin{eqnarray}\label{3.57}
 &&B_{2132}\nonumber\\
 &=&2r\sum_{k} g_n({\frac{k-1}{2^n}})\sum_{i=1}^m (\sigma^*\nabla\phi)_i(X({\frac{k-1}{2^n}})) (\sigma^*(X^n({\frac{k-1}{2^n}}))(X^n({\frac{k-1}{2^n}})-X({\frac{k-1}{2^n}}))_i\nonumber\\
 &&\times (|W_i({\frac{k}{2^n}})-W_i({\frac{k-1}{2^n}})|^2-\frac{1}{2^n})\nonumber\\
 &+&2r\sum_{k} g_n({\frac{k-1}{2^n}})\sum_{i\not= j}(\sigma^*\nabla\phi)_i(X({\frac{k-1}{2^n}})) (\sigma^*(X^n({\frac{k-1}{2^n}}))(X^n({\frac{k-1}{2^n}})-X({\frac{k-1}{2^n}}))_j\nonumber\\
 &&\times (W_i({\frac{k}{2^n}})-W_i({\frac{k-1}{2^n}}))(W_j({\frac{k}{2^n}})-W_j({\frac{k-1}{2^n}}))\nonumber\\
 &+&2r\sum_{k} g_n({\frac{k-1}{2^n}})\sum_{i=1}^m (\sigma^*\nabla\phi)_i(X({\frac{k-1}{2^n}})) (\sigma^*(X^n({\frac{k-1}{2^n}}))(X^n({\frac{k-1}{2^n}})-X({\frac{k-1}{2^n}}))_i(\frac{1}{2^n})\nonumber\\
 &:=& B_{21321}+B_{21322}+B_{21323}.
\end{eqnarray}
Conditioning on ${\cal F}_{\frac{k-1}{2^n}}$ and using the independence of $W_i$, $W_j$ for $i\not = j$, we see that $E[B_{21321}]=0$ and  $E[B_{21322}]=0$.
Furthermore,
\begin{eqnarray}\label{3.58}
 &&E[B_{21323}]\nonumber\\
 &=&2rE[\int_0^tg_n(s)<\sigma^*\nabla\phi(X(s)),  \sigma^*(X^n(s))(X^n(s)-X(s))>ds]\nonumber\\
 &+&2rE[\int_0^t\{g_n(s_n^-)<\sigma^*\nabla\phi(X(s_n^-)),  \sigma^*(X^n(s_n^-))(X^n(s_n^{-})-X(s_n^{-}))>\nonumber\\
 && -g_n(s)<\sigma^*\nabla\phi(X(s)),  \sigma^*(X^n(s))(X^n(s)-X(s))>\} ds]\nonumber\\
 &\leq& 2rE[\int_0^tg_n(s)<\sigma^*\nabla\phi(X(s)),  \sigma^*(X^n(s))(X^n(s)-X(s))>ds]\nonumber\\
 &&\quad\quad +C(\frac{1}{2^n})^{\frac{1}{2}},
\end{eqnarray}
where (\ref{3.01}) and (\ref{3.02}) were again used.
Combining together (\ref{3.53})---(\ref{3.58}) we obtain that
\begin{eqnarray}\label{3.59}
 &&E[B_{21}]\nonumber\\
 &\leq& 2rE[\int_0^tg_n(s)<\sigma^*\nabla\phi(X(s)),  \sigma^*(X^n(s))(X^n(s)-X(s))>ds]\nonumber\\
 &&\quad\quad +C(\frac{1}{2^n})^{\frac{1}{2}}.
\end{eqnarray}
To bound the term $B_{26}$, we write it as
\begin{eqnarray}\label{3.60}
 &&B_{26}\nonumber\\
 &=&
2r\int_0^t(\int_{s_n^-}^s[g_n(u)<\nabla\phi(X^n(u)), \sigma(X^n(u)dW^n(u)>-g_n(s_n^-)<\nabla\phi(X^n(s_n^-)),\nonumber\\
&&\quad\quad  \sigma(X^n(s_n^-))dW^n(u)>])<X^n(s)-X(s)
, \sigma(X^n(s))dW^n(s)> \nonumber\\
&+& 2r\int_0^t(\int_{s_n^-}^sg_n(s_n^-)<\nabla\phi(X^n(s_n^-)), \sigma(X^n(s_n^-))dW^n(u)>)\nonumber\\
&&\times [<\sigma^*(X^n(s))(X^n(s)-X(s))-\sigma^*(X^n(s_n^-))(X^n(s_n^-)-X(s_n^-)), dW^n(s)>]\nonumber\\
&+& 2r\int_0^tg_n(s_n^-)<\nabla\phi(X^n(s_n^-)), \sigma(X^n(s_n^-))(W^n(s)-W^n(s_n^-))>\nonumber\\
&&\quad\quad \times <\sigma^*(X^n(s_n^-))(X^n(s_n^-)-X(s_n^-)), dW^n(s)>\nonumber\\
&:=&B_{261}+B_{262}+B_{263}.
\end{eqnarray}
Following the same arguments leading to the estimates for $B_{211}$, $B_{212}$, it can be shown that
\begin{equation}\label{3.61}
E[B_{261}]\leq C(\frac{1}{2^n})^{\frac{1}{2}}, \quad\quad E[B_{261}]\leq C(\frac{1}{2^n})^{\frac{1}{2}}.
\end{equation}
and
\begin{eqnarray}\label{3.62}
 &&E[B_{263}]\nonumber\\
 &\leq& rE[\int_0^tg_n(s)<\sigma^*\nabla\phi(X^n(s)),  \sigma^*(X^n(s))(X^n(s)-X(s))>ds]\nonumber\\
 &&\quad\quad +C(\frac{1}{2^n})^{\frac{1}{2}}.
\end{eqnarray}
Putting together (\ref{3.48})----(\ref{3.62}) we get
\begin{eqnarray}\label{3.62-1}
 &&E[B_2]\nonumber\\
 &\leq& rE[\int_0^tg_n(s)<\sigma^*\nabla\phi(X^n(s)),  \sigma^*(X^n(s))(X^n(s)-X(s))>ds]\nonumber\\
 &&+2rE[\int_0^tg_n(s)<\sigma^*\nabla\phi(X(s)),  \sigma^*(X^n(s))(X^n(s)-X(s))>ds]\nonumber\\
 &&\quad\quad +C(\frac{1}{2^n})^{\frac{1}{2}}.
\end{eqnarray}

Now we turn to the $B_3$. Using the equations satisfied by $X^n$ and $X$, we have
\begin{eqnarray}\label{3.63}
 B_3
 &=&
2\int_0^tg_n(s_n^{-})<\int_{s_n^-}^s\sigma(X^n(u))dW^n(u),\sigma(X^n(s))dW^n(s)> \nonumber\\
&+& 2\int_0^tg_n(s_n^{-})<L^n(s)-L^n(s_n^{-}),\sigma(X^n(s))dW^n(s)>\nonumber\\
&+&2\int_0^tg_n(s_n^{-})<\int_{s_n^-}^sb(X^n(u))du,\sigma(X^n(s))dW^n(s)> \nonumber\\
&-&
2\int_0^tg_n(s_n^{-})<\int_{s_n^-}^s\sigma(X(u))dW(u),\sigma(X^n(s))dW^n(s)> \nonumber\\
&-&
\int_0^tg_n(s_n^{-})<\int_{s_n^-}^s\sigma\sigma^{\prime}(X(u))du,\sigma(X^n(s))dW^n(s)> \nonumber\\
&-&
2\int_0^tg_n(s_n^{-})<\int_{s_n^-}^sb(X(u))du,\sigma(X^n(s))dW^n(s)> \nonumber\\
&-& 2\int_0^tg_n(s_n^{-})<L(s)-L(s_n^{-}),\sigma(X^n(s))dW^n(s)>\nonumber\\
&:=&\sum_{i=1}^7B_{3i}.
\end{eqnarray}
Similar to the term $A_{25}$, we have
\begin{eqnarray}\label{3.64}
 &&E[B_{32}]\nonumber\\
 &\leq &C E[\sum_{k} \int_{\frac{k}{2^n}}^{\frac{k+1}{2^n}}g_n(\frac{k-1}{2^n})(\int_{\frac{k-1}{2^n}}^{s}\nu(X^n(u))d|L^n|(u))
 2^n|W(\frac{k}{2^n})-W(\frac{k-1}{2^n})|ds]\nonumber\\
 &\leq& 2CE[|L^n|(t) \sup_{|u-v|\leq \frac{1}{2^n}}(|W(u)-W(v)|)]\nonumber\\
 &\leq& 2C(E[|L^n|^2(t)])^{\frac{1}{2}}(\frac{1}{2^n})^{\frac{1}{2}}\leq C(\frac{1}{2^n})^{\frac{1}{2}}.
\end{eqnarray}
By the same reason,
\begin{equation}\label{3.65}
E[B_{37}]\leq C(\frac{1}{2^n})^{\frac{1}{2}}.
\end{equation}
Using a similar argument as for the term $A_{22}$, we obtain
\begin{equation}\label{3.66}
E[B_{3i}]\leq C(\frac{1}{2^n})^{\frac{1}{2}}, i=3, 5, 6.
\end{equation}
To bound  $B_{31}$, we write it as
\begin{eqnarray}\label{3.67}
 &&B_{31}\nonumber\\
 &=& 2\int_0^tg_n(s_n^{-})<\int_{s_n^-}^s(\sigma(X^n(u))-\sigma(X^n(s_n^{-})))dW^n(u),\sigma(X^n(s))dW^n(s)> \nonumber\\
 &+& 2\int_0^tg_n(s_n^{-})<\sigma(X^n(s_n^{-})))(W^n(s)-W^n(s_n^{-})),(\sigma(X^n(s))-\sigma(X^n(s_n^{-})))dW^n(s)> \nonumber\\
 &+& 2\int_0^tg_n(s_n^{-})<\sigma(X^n(s_n^{-}))(W^n(s)-W^n(s_n^{-})),\sigma(X^n(s_n^{-}))dW^n(s)> \nonumber\\
 &:=& B_{311}+B_{312}+B_{313},
 \end{eqnarray}
 where
 \begin{eqnarray}\label{3.68}
 &&E[B_{311}]\nonumber\\
 &\leq & C\int_0^tds\int_{s_n^-}^sE[|X^n(u)-X^n(s_n^{-})|\dot{W}^n(u)||\dot{W}^n(s)| ]du\nonumber\\
 &\leq &  C(\frac{1}{2^n})^{\frac{1}{2}},
 \end{eqnarray}
 \begin{eqnarray}\label{3.69}
 &&E[B_{312}]\nonumber\\
 &\leq & C\int_0^tdsE[|X^n(s)-X^n(s_n^{-})||W^n(s)-W^n(s_n^{-})||\dot{W}^n(s)|] \nonumber\\
 &\leq &  C(\frac{1}{2^n})^{\frac{1}{2}},
 \end{eqnarray}
 and
 \begin{eqnarray}\label{3.70}
 &&B_{313}\nonumber\\
 &=&2\sum_{k}2^{2n} \int_{\frac{k}{2^n}}^{\frac{k+1}{2^n}}g_n(\frac{k-1}{2^n})(s-\frac{k}{2^n})ds |\sigma(X^n(\frac{k-1}{2^n})) (W(\frac{k}{2^n})-W(\frac{k-1}{2^n}))|^2\nonumber\\
 &+& 2\sum_{k} g_n(\frac{k-1}{2^n})<\sigma(X^n(\frac{k-1}{2^n})) (W(\frac{k-1}{2^n})-W(\frac{k-2}{2^n})), \nonumber\\
 &&\quad\quad \sigma(X^n(\frac{k-1}{2^n}))(W(\frac{k}{2^n})-W(\frac{k-1}{2^n}))>\nonumber\\
 &:=&B_{3131}+B_{3132}.
\end{eqnarray}
Conditioning on ${\cal F}_{\frac{k-1}{2^n}}$, we see that $E[B_{3132}]=0$. Rearranging  the terms, we find that
\begin{eqnarray}\label{3.71}
 &&B_{3131}\nonumber\\
 &=&\sum_{k}g_n(\frac{k-1}{2^n})\sum_{i=1}^d\sum_{j=1}^m\sigma_{ij}^2(X^n(\frac{k-1}{2^n})\{ (W_j(\frac{k}{2^n})-W_j(\frac{k-1}{2^n}))^2-\frac{1}{2^n}\}\nonumber\\
 &+&\sum_{k}g_n(\frac{k-1}{2^n})\sum_{i=1}^d\sum_{j\not =l}^m\sigma_{ij}(X^n(\frac{k-1}{2^n}))\sigma_{il}(X^n(\frac{k-1}{2^n}) \nonumber\\
 &&\quad\times (W_j(\frac{k}{2^n})-W_j(\frac{k-1}{2^n}))(W_l(\frac{k}{2^n})-W_l(\frac{k-1}{2^n}))\nonumber\\
 &+&\int_0^t\{g_n(s_n^-)\sum_{i=1}^d\sum_{j=1}^m\sigma_{ij}^2(X^n(s_n^-))-g_n(s)\sum_{i=1}^d\sum_{j=1}^m\sigma_{ij}^2(X^n(s))\}ds\nonumber\\
 &+&\int_0^tg_n(s)\sum_{i=1}^d\sum_{j=1}^m\sigma_{ij}^2(X^n(s))ds.
 \end{eqnarray}
 By conditioning and using the independence of $W_j$ and $W_l$ for $j\not =l$, we see that the expectation of the first two terms on the right side are zero. The expectation of the third term is bounded by $C(\frac{1}{2^n})^{\frac{1}{2}}$. Thus we have
 \begin{eqnarray}\label{3.72}
 &&E[B_{3131}]\nonumber\\
 &\leq &E[\int_0^tg_n(s)\sum_{i=1}^d\sum_{j=1}^m\sigma_{ij}^2(X^n(s))ds] +C(\frac{1}{2^n})^{\frac{1}{2}}.
 \end{eqnarray}
 As for  $B_{34}$, we have
\begin{eqnarray}\label{3.67-1}
 &&B_{34}\nonumber\\
 &=& -2\int_0^tg_n(s_n^{-})<\int_{s_n^-}^s(\sigma(X(u))-\sigma(X(s_n^{-})))dW(u),\sigma(X^n(s))dW^n(s)> \nonumber\\
 &-& 2\int_0^tg_n(s_n^{-})<\sigma(X(s_n^{-})))(W(s)-W(s_n^{-})),(\sigma(X^n(s))-\sigma(X^n(s_n^{-})))dW^n(s)> \nonumber\\
 &-& 2\int_0^tg_n(s_n^{-})<\sigma(X(s_n^{-}))(W(s)-W(s_n^{-})),\sigma(X^n(s_n^{-}))dW^n(s)> \nonumber\\
 &:=& B_{341}+B_{342}+B_{343},
 \end{eqnarray}
 Similar to the terms $B_{311},B_{312}$, we have
 \begin{eqnarray}\label{3.68-1}
 E[B_{341}]
 &\leq &  C(\frac{1}{2^n})^{\frac{1}{2}},
 \end{eqnarray}
 \begin{eqnarray}\label{3.69-1}
 E[B_{342}]
 &\leq &  C(\frac{1}{2^n})^{\frac{1}{2}}.
 \end{eqnarray}
 Moreover,
 \begin{eqnarray}\label{3.70-1}
 &&B_{343}\nonumber\\
 &=&-2\sum_{k} g_n(\frac{k-1}{2^n})<\sigma(X(\frac{k-1}{2^n})) (W(\frac{k}{2^n})-W(\frac{k-1}{2^n})),\nonumber\\
 &&\quad \quad \sigma(X^n(\frac{k-1}{2^n})) (W(\frac{k}{2^n})-W(\frac{k-1}{2^n}))>\nonumber\\
 &-& 2\sum_{k}2^n \int_{\frac{k}{2^n}}^{\frac{k+1}{2^n}}ds g_n(\frac{k-1}{2^n})<\sigma(X(\frac{k-1}{2^n}) (W(s)-W(\frac{k}{2^n})), \nonumber\\
 &&\quad\quad \sigma(X^n(\frac{k-1}{2^n}) (W(\frac{k}{2^n})-W(\frac{k-1}{2^n}))>\nonumber\\
 &:=&B_{3431}+B_{3432}.
\end{eqnarray}
Conditioning on ${\cal F}_{\frac{k}{2^n}}$, we see that $E[B_{3432}]=0$. Term $B_{3431}$ can be written as follows:
\begin{eqnarray}\label{3.71-1}
 &&B_{3431}\nonumber\\
 &=&-2\sum_{k}g_n(\frac{k-1}{2^n})\sum_{i=1}^d\sum_{j=1}^m\sigma_{ij}^2(X(\frac{k-1}{2^n}))\sigma_{ij}(X^n(\frac{k-1}{2^n}))\nonumber\\
 &&\quad\quad \{ (W_j(\frac{k}{2^n})-W_j(\frac{k-1}{2^n}))^2-\frac{1}{2^n}\}\nonumber\\
 &-2&\sum_{k}g_n(\frac{k-1}{2^n})\sum_{i=1}^d\sum_{j\not =l}^m\sigma_{ij}(X(\frac{k-1}{2^n}))\sigma_{il}(X^n(\frac{k-1}{2^n})) \nonumber\\
 &&\quad\times (W_j(\frac{k}{2^n})-W_j(\frac{k-1}{2^n}))(W_l(\frac{k}{2^n})-W_l(\frac{k-1}{2^n}))\nonumber\\
 &-2&\int_0^t\{g_n(s_n^-)\sum_{i=1}^d\sum_{j=1}^m\sigma_{ij}^2(X(s_n^-))\sigma_{ij}^2(X^n(s_n^-))-g_n(s)\sum_{i=1}^d\sum_{j=1}^m
 \sigma_{ij}^2(X(s))\sigma_{ij}^2(X^n(s))\}ds\nonumber\\
 &-2&\int_0^tg_n(s)\sum_{i=1}^d\sum_{j=1}^m
 \sigma_{ij}(X(s))\sigma_{ij}(X^n(s))ds.
 \end{eqnarray}
 By conditioning and using the independence of $W_j$ and $W_l$ for $j\not =l$, it follows that the expectation of the first two terms on the right side are zero and the expectation of the third term is bounded by $C(\frac{1}{2^n})^{\frac{1}{2}}$. Thus we have
 \begin{eqnarray}\label{3.72-2}
 &&E[B_{3431}]\nonumber\\
 &\leq &E[\int_0^tg_n(s)\sum_{i=1}^d\sum_{j=1}^m\sigma_{ij}^2(X^n(s))ds] +C(\frac{1}{2^n})^{\frac{1}{2}}.
 \end{eqnarray}
 It follows from (\ref{3.63}) ---(\ref{3.72-2}) that
 \begin{eqnarray}\label{3.72-1}
 &&E[B_{3}]\nonumber\\
 &\leq &E[\int_0^tg_n(s)\sum_{i=1}^d\sum_{j=1}^m\sigma_{ij}^2(X^n(s))ds] +C(\frac{1}{2^n})^{\frac{1}{2}}\nonumber\\
 &-&2E[\int_0^tg_n(s)\sum_{i=1}^d\sum_{j=1}^m
 \sigma_{ij}(X(s))\sigma_{ij}(X^n(s))ds].
 \end{eqnarray}
 To bound  $B_4$, denote by $\nabla\sigma =(\nabla \sigma_{ij})\in R^{d\times m}\otimes R^d$ and  $\sigma^*\nabla\sigma \in R^{d\times m}\otimes R^d$ the linear mappings defined by
 $$<\nabla \sigma, y>=(<\nabla\sigma_{ij},y>)\in R^{d\times m}, \quad y\in R^d,$$
 $$<\sigma^*\nabla \sigma, x>=(<\sigma^*\nabla\sigma_{ij},x>)\in R^{d\times m}, \quad x\in R^m.$$
 Observe that
 \begin{eqnarray}\label{3.73}
B_4&=& 2 \int_0^tg_n(s_n^{-})<X^n(s_n^{-})-X(s_n^{-}),(\int_{s_n^-}^s<\nabla\sigma(X^n(u)),\sigma(X^n(u))dW^n(u)>)dW^n(s)> \nonumber\\
&+& 2 \int_0^tg_n(s_n^{-})<X^n(s_n^{-})-X(s_n^{-}),(\int_{s_n^-}^s<\nabla\sigma(X^n(u)),\nu(X^n(u))d|L^n|_u>)dW^n(s)> \nonumber\\
&+& 2 \int_0^tg_n(s_n^{-})<X^n(s_n^{-})-X(s_n^{-}),(\int_{s_n^-}^s<\nabla\sigma(X^n(u)),b(X^n(u))du>)dW^n(s)> \nonumber\\
&:=& B_{41}+B_{42}+B_{43}.
\end{eqnarray}
As other similar terms treated above, we can show that
\begin{equation}\label{3.74}
E[B_{43}]\leq C(\frac{1}{2^n})^{\frac{1}{2}}.
\end{equation}
\begin{equation}\label{3.75}
E[B_{42}]\leq CE[|L^n|_t \sup_{|u-v|\leq \frac{1}{2^n}}|W(u)-W(v)|]\leq C(\frac{1}{2^n})^{\frac{1}{2}}.
\end{equation}
$B_{41}$ can be further split as
\begin{eqnarray}\label{3.76}
B_{41}&=& 2 \int_0^tg_n(s_n^{-})<X^n(s_n^{-})-X(s_n^{-}),\nonumber\\
&&\quad\quad (\int_{s_n^-}^s<\sigma^*\nabla\sigma(X^n(u))-\sigma^*\nabla\sigma(X^n(s_n^-)),dW^n(u)>)dW^n(s)> \nonumber\\
&+&2\sum_{k}g_n(\frac{k-1}{2^n})<X^n(\frac{k-1}{2^n})-X(\frac{k-1}{2^n}), \nonumber\\
&&<\sigma^*\nabla\sigma(X^n(\frac{k-1}{2^n})), W(\frac{k-1}{2^n})-W(\frac{k-2}{2^n})>(W(\frac{k}{2^n})-W(\frac{k-1}{2^n}))>\nonumber\\
&+&2\sum_{k}2^{2n} \int_{\frac{k}{2^n}}^{\frac{k+1}{2^n}}(s-{\frac{k}{2^n}})ds g_n(\frac{k-1}{2^n})<X^n(\frac{k-1}{2^n})-X(\frac{k-1}{2^n}), \nonumber\\
&&<\sigma^*\nabla\sigma(X^n(\frac{k-1}{2^n})), W(\frac{k}{2^n})-W(\frac{k-1}{2^n})>(W(\frac{k}{2^n})-W(\frac{k-1}{2^n}))>\nonumber\\
&:=& B_{411}+B_{412}+B_{413}.
\end{eqnarray}
As other similar terms above, we have
\begin{equation}\label{3.77}
E[B_{411}]\leq C(\frac{1}{2^n})^{\frac{1}{2}}, \quad \quad E[B_{412}]=0.
\end{equation}
On the other hand, we have
\begin{eqnarray}\label{3.78}
B_{413}&=& \sum_{k} g_n(\frac{k-1}{2^n})\sum_{i=1}^d(X_i^n(\frac{k-1}{2^n})-X_i(\frac{k-1}{2^n}))\nonumber\\
&&  \quad \times \sum_{j\not =l}^m (\sigma^*\nabla\sigma_{ij})_l(X^n(\frac{k-1}{2^n})) (W_l(\frac{k}{2^n})-W_l(\frac{k-1}{2^n}))(W_j(\frac{k}{2^n})-W_j(\frac{k-1}{2^n}))\nonumber\\
&+&
\sum_{k} g_n(\frac{k-1}{2^n})\sum_{i=1}^d(X_i^n(\frac{k-1}{2^n})-X_i(\frac{k-1}{2^n}))\nonumber\\
&&  \quad \times \sum_{j=1}^m (\sigma^*\nabla\sigma_{ij})_j(X^n(\frac{k-1}{2^n})) \{|W_j(\frac{k}{2^n})-W_j(\frac{k-1}{2^n})|^2-\frac{1}{2^n}\}\nonumber\\
&+& \int_0^t\{ g_n(s_n^-)\sum_{i=1}^d(X_i^n(s_n^-)-X_i(s_n^-))\sum_{j =1}^m (\sigma^*\nabla\sigma_{ij})_j(X^n(s_n^-))\nonumber\\
&& \quad - g_n(s)\sum_{i=1}^d(X_i^n(s)-X_i(s))\sum_{j =1}^m (\sigma^*\nabla\sigma_{ij})_j(X^n(s))\}ds\nonumber\\
&+& \int_0^tg_n(s)\sum_{i=1}^d(X_i^n(s)-X_i(s))\sum_{j =1}^m (\sigma^*\nabla\sigma_{ij})_j(X^n(s))ds.
\end{eqnarray}
Using the independence of $W_j$ and $W_l$ for $j\not =l$, (\ref{3.01}) and (\ref{3.02}), by conditioning on ${\cal F}_{\frac{k-1}{2^n}}$ we obtain from
(\ref{3.78}) that
\begin{eqnarray}\label{3.79}
&&E[B_{413}]\nonumber\\
&\leq & E[\int_0^tg_n(s)\sum_{i=1}^d(X_i^n(s)-X_i(s))\sum_{j =1}^m (\sigma^*\nabla\sigma_{ij})_j(X^n(s))\}ds]\nonumber\\
&&+C(\frac{1}{2^n})^{\frac{1}{2}}.
\end{eqnarray}
Combining (\ref{3.73})---(\ref{3.79}) yields that
\begin{eqnarray}\label{3.80}
&&E[B_{4}]\nonumber\\
&\leq & E[\int_0^tg_n(s)\sum_{i=1}^d(X_i^n(s)-X_i(s))\sum_{j =1}^m (\sigma^*\nabla\sigma_{ij})_j(X^n(s))\}ds]\nonumber\\
&&+C(\frac{1}{2^n})^{\frac{1}{2}}.
\end{eqnarray}
The lemma is proved by putting together (\ref{3.47}), (\ref{3.62-1}), (\ref{3.72-1}) and  (\ref{3.82}).
\vskip 0.3cm
\noindent{\bf Proof of Theorem 2.2: (Continued)}.
Choose $r<-\frac{2C_0}{\alpha}$, where $\alpha$, $C_0$ are the constants appeared in the assumptions (D.1) and (D.2).
By the Lipschitz continuity of the coefficients and boundedness of $\phi$, $\phi^{\prime\prime}$, $\nabla \phi$, $\sigma\sigma^{\prime}$ on
the domain $\bar{D}$, it follows from (\ref{3.1}) that
\begin{eqnarray}\label{3.81}
&& E[f_n(t)]\nonumber\\
&\leq & C_rE[\int_0^tf_n(s)ds]\nonumber\\
&+& E[\int_0^t\{<rf_n(s)\nabla\phi(X(s))-2g_n(s)(X^n(s)-X(s)), \nu(X(s))>\}d|L|(s)]\nonumber\\
&+&rE[\int_0^tf_n(s)<\nabla\phi(X^n(s)), \sigma(X^n(s))dW^n(s)>]\nonumber\\
&+&E[\int_0^t\{<rf_n(s)\nabla\phi(X^n(s))+2g_n(s)(X^n(s)-X(s)), \nu(X^n(s))>\}d|L^n|(s)]\nonumber\\
&+&2E[\int_0^tg_n(s)<X^n(s)-X(s), \sigma(X^n(s))dW^n(s)>]\nonumber\\
&-&E[\int_0^tg_n(s)<X^n(s)-X(s), \sigma\sigma^{\prime}(X(s))>ds]\nonumber\\
&+&E[\int_0^tg_n(s) tr(\sigma\sigma^*(X(s)))ds]\nonumber\\
&-&2rE[\int_0^tg_n(s)<\sigma^*(X(s))(X^n(s)-X(s)), \sigma^*\nabla\phi(X(s))>ds].
\end{eqnarray}
In view of $r<0$ and the assumptions (D.1) and (D.2), we deduce that
\begin{eqnarray}\label{3.82}
&& <rf_n(s)\nabla\phi(X(s))-2g_n(s)(X^n(s)-X(s)), \nu(X(s))>\nonumber\\
&=&g_n(s)[r<\nabla\phi(X(s)),\nu(X(s))>|x^n(s)-X(s)|^2-2<X^n(s)-X(s), \nu(X(s))>]\nonumber\\
&\leq &g_n(s)[r\alpha |x^n(s)-X(s)|^2+2C_0|X^n(s)-X(s)|^2]\leq 0,
\end{eqnarray}
and similarly
\begin{eqnarray}\label{3.83}
&& <rf_n(s)\nabla\phi(X^n(s))+2g_n(s)(X^n(s)-X(s)), \nu(X^n(s))>\nonumber\\
&\leq & 0.
\end{eqnarray}
Thus, using  Lemma 3.1 and Lemma 3.2, taking into account (\ref{3.82}) and (\ref{3.83}) we obtain from
(\ref{3.81}) that
\begin{eqnarray}\label{3.84}
&& E[f_n(t)]\nonumber\\
&\leq & C_rE[\int_0^tf_n(s)ds]+C(\frac{1}{2^n})^{\frac{1}{2}}\nonumber\\
&-&E[\int_0^tg_n(s)<X^n(s)-X(s), \sigma\sigma^{\prime}(X(s))>ds]\nonumber\\
&+&E[\int_0^tg_n(s) \sum_{i=1}^d\sum_{j=1}^m\sigma_{ij}^2(X(s))ds]\nonumber\\
&-&2rE[\int_0^tg_n(s)<\sigma^*(X(s))(X^n(s)-X(s)), \sigma^*\nabla\phi(X(s))>ds]\nonumber\\
&+& 2r\int_0^t <g_n(s)\sigma^*(X^n(s))(X^n(s)
 -X(s)), \sigma^*\nabla\phi(X^n(s))>ds\nonumber\\
 &-& 2r\int_0^t <g_n(s)\sigma^*(X(s))(X^n(s)
 -X(s)), \sigma^*\nabla\phi(X^n(s))>ds\nonumber\\
 &+&2rE[\int_0^tg_n(s)<\sigma^*\nabla\phi(X(s)),  \sigma^*(X^n(s))(X^n(s)-X(s))>ds]\nonumber\\
 &+&E[\int_0^tg_n(s)\sum_{i=1}^d\sum_{j=1}^m\sigma_{ij}^2(X^n(s))ds]\nonumber\\
 &+&E[\int_0^tg_n(s)<X^n(s)-X(s), \sigma\sigma^{\prime}(X^n(s))>ds]\nonumber\\
 \nonumber\\
 &-2&E[\int_0^tg_n(s)\sum_{i=1}^d\sum_{j=1}^m\sigma_{ij}(X(s)\sigma_{ij}(X^n(s))ds]\nonumber\\
 &\leq & C_rE[\int_0^tf_n(s)ds]+C(\frac{1}{2^n})^{\frac{1}{2}}\nonumber\\
&+&E[\int_0^tg_n(s) \sum_{i=1}^d\sum_{j=1}^m(\sigma_{ij}(X(s))-\sigma_{ij}(X^n(s)))^2ds]\nonumber\\
&+& 2r\int_0^t <g_n(s)(\sigma^*(X^n(s))-\sigma^*(X(s)))(X^n(s)
 -X(s)),\nonumber\\
 &&\quad\quad\quad  \sigma^*\nabla\phi(X^n(s))>ds\nonumber\\
 &+&2rE[\int_0^tg_n(s)<\sigma^*\nabla\phi(X(s)), \nonumber\\
 &&\quad\quad (\sigma^*(X^n(s))-\sigma^*(X(s)))(X^n(s)-X(s))>ds]\nonumber\\
 &+&E[\int_0^tg_n(s)<X^n(s)-X(s), \sigma\sigma^{\prime}(X^n(s))-\sigma\sigma^{\prime}(X(s))>ds]\nonumber\\
 &\leq&CE[\int_0^tf_n(s)ds]+C(\frac{1}{2^n})^{\frac{1}{2}},
 \end{eqnarray}
 where the Lipschitz continuity of the coefficients and the fact that $f_n(s)=g_n(s)|X^n(s)-X(s)|^2$ have been used. Finally by the Gronwall's inequality, we obtain
 \begin{equation}\label{3.85}
 E[f_n(t)]\leq C(\frac{1}{2^n})^{\frac{1}{2}}\rightarrow 0
 \end{equation}
 as $n\rightarrow \infty$, completing the proof of (\ref{3.05}), hence  the theorem.

\end{document}